%% file: main.tex
\documentclass[a4paper,11pt]{article}
\usepackage[utf8]{inputenc}
\usepackage{a4wide}

\usepackage{amsmath}
\usepackage{amssymb}
\usepackage{amsthm}
\usepackage{graphicx}
\usepackage{mathtools}
\usepackage{url}
\usepackage{multirow}
\usepackage{xcolor}
\usepackage{caption}
\usepackage{subcaption}

\definecolor{darkblue}{RGB}{0,60,180}
\definecolor{darkgreen}{RGB}{0,130,70}
\definecolor{darkorange}{RGB}{180,60,0}

\providecommand{\keywords}[1]
{
  \noindent \small
  \textbf{Keywords:} #1
}
\providecommand{\amscode}[1]
{
  \noindent \small
  \textbf{AMS subject classifications:} #1
}

\mathtoolsset{showonlyrefs}
\allowdisplaybreaks[1]

\newtheorem{problem}{Problem}
\newtheorem{theorem}{Theorem}

\newtheorem{algorithm}{Algorithm}
\newtheorem{lemma}{Lemma}
\newtheorem{remark}{Remark}

\title{Estimating nested expectations without inner conditional sampling and application to value of information analysis}
\author{Tomohiko Hironaka\thanks{School of Engineering, University of Tokyo, 7-3-1 Hongo, Bunkyo-ku, Tokyo 113-8656, Japan ({\tt hironaka-tomohiko@g.ecc.u-tokyo.ac.jp}; {\tt goda@frcer.t.u-tokyo.ac.jp})}, Takashi Goda\footnotemark[1]}

\begin{document}

\maketitle
\begin{abstract}
    Motivated by various computational applications, we investigate the problem of estimating nested expectations. Building upon recent work by the authors, we propose a novel Monte Carlo estimator for nested expectations, inspired by sparse grid quadrature, that does not require sampling from inner conditional distributions. Theoretical analysis establishes an upper bound on the mean squared error of our estimator under mild assumptions on the problem, demonstrating its efficiency for cases with low-dimensional outer variables. We illustrate the effectiveness of our estimator through its application to problems related to value of information analysis, with moderate dimensionality. Overall, our method presents a promising approach to efficiently estimate nested expectations in practical computational settings.
\end{abstract}

\keywords{Monte Carlo methods, nested expectations, sparse grid quadrature, value of information analysis, mean squared error}

\amscode{65C05, 65D32, 65D40, 90B50}

\input{10-introduction}

\input{20-method}

\input{30-theoretical_th}

\input{40-numerical}
\input{60-conclusion}
\input{90-backmatter}

\end{document}

%% file: 10-introduction.tex
\section{Introduction}\label{sec:intro}

\subsection{Overview}
Estimating nested expectations is a computationally demanding problem that arises in various computational science and engineering fields. In many applications, it is necessary to compute the expectation of a function with respect to random variables, where the argument is given by the inner conditional expectation of other random variables. This type of quantity is called a \emph{nested expectation}, and it arises in a variety of fields, including Bayesian experimental designs \cite{lindley1956measure,chaloner1995bayesian}, risk estimation in computational finance \cite{glasserman2004monte,gordy2010nested,broadie2011efficient}, global sensitivity analysis \cite{sobol2001global,saltelli2008global}, and value of information analysis \cite{raiffa1961applied,howard1966information,welton2012evidence}. The standard Monte Carlo approach to estimate nested expectations simply nests the Monte Carlo estimators for inner and outer expectations, but it requires sampling from inner conditional distributions. This can be difficult or even impossible, and even if this is possible, the nested Monte Carlo estimator is known to be inefficient in terms of error convergence \cite{ryan2003estimating,ades2004expected,rainforth2018nesting}. To overcome this challenge, various methods have been proposed, including importance sampling, Markov chain Monte Carlo, and multilevel Monte Carlo methods. See \cite{Long2013Fast,ryan2016review,beck2018fast,beck2019multilevel,giles2019decision,giles2019multilevel,goda2020multilevel,hironaka2020multilevel,goda2022unbiased,wang2022unbiased} among many others. Nonetheless, most of these methods still rely on inner conditional sampling.

In this paper, we propose a novel Monte Carlo estimator for nested expectations that does not require sampling from inner conditional distributions. Our method, introduced in Section~\ref{sec:method}, is built upon recent work by the authors \cite{hironaka2023efficient} and is inspired by sparse grid quadrature, which is a numerical method for approximating high-dimensional integrals \cite{smolyak1963quadrature,Gerstner1998Numerical,bungartz2004sparse}. While our proposed estimator is theoretically shown to be efficient for problems with low-dimensional outer variables, we also empirically demonstrate its applicability to problems with moderate dimensionality. Specifically, in Section~\ref{sec:theoretical}, we provide a theoretical analysis that establishes an upper bound on the mean squared error (MSE) of our estimator under mild assumptions on the problem, and in Section~\ref{sec:numerical}, we demonstrate its effectiveness through numerical experiments related to value of information analysis. Our proposed method presents a promising approach to efficiently estimate nested expectations in practical computational settings.

\subsection{Related Works}
Before moving on, we would like to highlight some relevant works on estimating nested expectations without relying on inner conditional sampling. In the realm of global sensitivity analysis, Broto et al.\ proposed one approach that uses the $k$-nearest neighbors algorithm ($k$NN) to estimate sensitivity indices \cite{broto2020variance}. In the context of value of information analysis for medical decision problems, several papers have proposed different approaches to estimate the expected value of partial perfect information (EVPPI) and the expected value of sample information (EVSI), including \cite{strong2014estimating,strong2015estimating,menzies2016efficient,heath2018efficient,jalal2018gaussian,heath2020calculating}. For instance, Strong et al.\ suggested a regression-based approach using generalized additive models or Gaussian process regression in \cite{strong2014estimating,strong2015estimating}. The combined use of Monte Carlo sampling and regression has also been studied independently in the context of financial engineering, as in \cite{longstaff2001valuing,broadie2015risk}. Furthermore, Hong et al.\ investigated the kernel smoothing approach as well as the $k$NN approach in the context of risk estimation in financial engineering \cite{hong2017kernel}. The key, common point among these studies as well as this paper is that nested expectations are estimated through post-computation using only random samples from the joint distribution of the inner and outer random variables, without explicitly sampling from the inner conditional distributions.

In a recent work by the authors \cite{hironaka2023efficient}, we have introduced a Monte Carlo approach different from the aforementioned works. Given a set of random samples from the joint distribution of the inner and outer random variables, our approach stratifies them a posteriori into mutually disjoint subsets, depending on the values of the outer random variables. For each subset, the inner conditional expectation is estimated by taking the average of the sample values of the inner random variables that belong to that subset. Although we established an upper bound on the MSE, we made a relatively strong assumption that the outer random variables are all continuous and mutually independent. Additionally, our estimator requires the number of random samples to be of the form $m^{2K}$, where $K$ is the number of outer random variables and $m > 1$ is an integer. Therefore, an exponentially large number of random samples is needed as the outer variables' dimensionality increases. In contrast, the proposed method in this paper is free from these disadvantages: we establish an upper bound on the MSE, without assuming that the outer random variables are all continuous and mutually independent, for the number of random samples being any power of 2 irrespective of the dimension $K$.

\subsection{Problem Setting}
Let us now define nested expectations more precisely. Suppose we have two random variables $X\in \mathbb{R}^J$ and $Y\in \mathbb{R}^K$. For a function $f: \mathbb{R}^J \to \mathbb{R}$, the nested expectation is defined as
\begin{align}\label{eq:ne2}
I=\mathbb{E}_{\rho(Y)}f\left(\mathbb{E}_{\rho(X|Y)}X\right),
\end{align}
where $\rho(Y)$ denotes the marginal probability distribution of $Y$ and $\rho(X|Y)$ denotes the conditional probability distribution of $X$ given $Y$. Throughout this paper, we assume that it is possible to obtain independent and identically distributed (i.i.d.)\ samples from the joint distribution $\rho(X,Y)$ of the random variables $X$ and $Y$, although this may not be the case for the inner conditional distribution $\rho(X|Y)$. We further assume that evaluating the function $f$ pointwise and generating one sample from the joint distribution $\rho(X,Y)$ both take a unit computational cost. The problem addressed in this paper is to estimate $I$ efficiently under these assumptions.

As we have already explained, there are numerous applications that require computation of nested expectations. In the following, we present an example from the value of information analysis for medical decision problems, which is a particular application of interest in this paper.

Let $D$ be a finite set of potential medical treatments, each with uncertain outcomes and costs. We model the net benefit of treatment $d\in D$ as a function of the input random variables $\theta$, denoted by $\mathrm{NB}_d$. These variables include factors such as the probability of side effects and the cost of treatment. Our goal is to determine whether conducting a clinical trial or medical research to reduce the uncertainty of $\theta$ is worth the investment \cite{welton2012evidence}. EVSI is a measure of the cost-effectiveness of conducting such research. It calculates the average gain in net benefit from making the observation $Y$:
\begin{align}\label{eq:evsi}
\mathbb{E}_{\rho(Y)}\max_{d\in D}\mathbb{E}_{\rho(\theta|Y)}\mathrm{NB}_d(\theta)-\max_{d\in D}\mathbb{E}_{\rho(\theta)}\mathrm{NB}_d(\theta).
\end{align}
Here, the first term represents the average net benefit when choosing the optimal treatment based on the observation $Y$, and the second term represents the net benefit without making the observation. Notice that the first term is a nested expectation, as given in \eqref{eq:ne2}.

%% file: 20-method.tex
\section{Proposed Method}\label{sec:method}

Throughout this paper, we denote the set of positive integers by $\mathbb{N}$. For the sake of notation, we introduce $\mathbf{u}{:}p$ to denote the vector obtained by appending a real number $p$ to the end of a vector $\mathbf{u} \in \mathbb{R}^d$. Specifically, $\mathbf{u}{:}p = (u_1, \dots, u_d, p) \in \mathbb{R}^{d+1}$. With this notation in mind, we now introduce our new Monte Carlo estimator for nested expectations, which is an extension of the previous method by the authors \cite{hironaka2023efficient} and is inspired by sparse grid quadrature.

\begin{algorithm}\label{alg:proposed}
    Let $m \in \mathbb{N}$, and let $\{(X^{(i)}, Y^{(i)})\}_{1 \leq i \leq 2^m}$ be a set of i.i.d.\ random samples from the joint distribution $\rho(X,Y)$. For $d \in \{0, \dots, m\}$ and $\mathbf{u} \in \{0,1\}^d$, we construct $B^{d}_{\mathbf{u}}, B'^d_{\mathbf{u}} \subseteq \{1, \dots, 2^m\}$ recursively as follows: For $d=0$, let
    \begin{align}
        B^{0}_{\emptyset} = \{1, \dots, 2^m\}.
    \end{align}
    For $0 < d \leq m$, we separate $B^{d-1}_{\mathbf{u}}$ into two mutually disjoint sets in two different ways for each $\mathbf{u}\in \{0,1\}^{d-1}$ such that
    \begin{align}
    \quad
            B^{d}_{\mathbf{u}{:}0}
        \cup
            B^{d}_{\mathbf{u}{:}1}
    &=
            B'^{d}_{\mathbf{u}{:}0}
        \cup
            B'^{d}_{\mathbf{u}{:}1}
    =
        B^{d-1}_{\mathbf{u}},\\
    \quad
            B^{d}_{\mathbf{u}{:}0}
        \cap
            B^{d}_{\mathbf{u}{:}1}
    &=
            B'^{d}_{\mathbf{u}{:}0}
        \cap
            B'^{d}_{\mathbf{u}{:}1}
    =
        \emptyset, \\
        \left| B^{d}_{\mathbf{u}{:}0} \right|
    =
        \left| B^{d}_{\mathbf{u}{:}1} \right|
    &=
        \left| B'^{d}_{\mathbf{u}{:}0} \right|
    =
        \left| B'^{d}_{\mathbf{u}{:}1} \right|
    =
        2^{m-d},
    \end{align}
    where $Y^{(i_0)}_{(d-1 \bmod K)+1} \leq Y^{(i_1)}_{(d-1 \bmod K)+1}$
    for all $
        i_0 \in B^{d}_{\mathbf{u}{:}0}
    $ and  $
        i_1 \in B^{d}_{\mathbf{u}{:}1}
    $ and $i_0 \leq i_1$ for all $
        i_0 \in B'^{d}_{\mathbf{u}{:}0}
    $ and $
        i_1 \in B'^{d}_{\mathbf{u}{:}1}
    $. Then we estimate $I$ by
    \begin{align}\label{eq:proposed_estimator}
        \hat{I}
    &=
        \sum_{d=0}^{m}
        \frac{1}{2^d}
        \sum_{\mathbf{u} \in \{0, 1\}^d}
        f\left(
            \frac{1}{2^{m-d}}
            \sum_{i \in B^d_{\mathbf{u}}} X^{(i)}
        \right)
    -
        \sum_{d=1}^{m}
        \frac{1}{2^d}
        \sum_{\mathbf{u} \in \{0, 1\}^d}
        f\left(
            \frac{1}{2^{m-d}}
            \sum_{i \in B'^d_{\mathbf{u}}} X^{(i)}
        \right).
    \end{align}
\end{algorithm}

Here, we explain how this estimator is constructed. First, let us take a look at the first term of \eqref{eq:proposed_estimator}, which is a sum over $d$ from 0 to $m$. For the sake of notation, we denote each summand by
\[ \hat{P}^{(d,m-d)}:=\frac{1}{2^d}\sum_{\mathbf{u} \in \{0, 1\}^d}f\left(\frac{1}{2^{m-d}}\sum_{i \in B^d_{\mathbf{u}}} X^{(i)}\right),\]
in which, the outer expectation is approximated by the average of $2^d$ samples over $\mathbf{u}$, and the inner conditional expectation is approximated by the average of $2^{m-d}$ samples in $B_{\mathbf{u}}^d$ for the corresponding $\mathbf{u}$. Here we note that $B_{\mathbf{u}}^d$ is defined by stratifying the samples $\{(X^{(i)}, Y^{(i)})\}_{1 \leq i \leq 2^m}$ depending on the values of $Y^{(i)}$, so that each sample in the same $B_{\mathbf{u}}^d$ is ``close'' to each other, which makes it possible to estimate the inner conditional expectation by the average over $B_{\mathbf{u}}^d$ without explicitly relying on sampling from the inner conditional distribution. Thus, as $d$ increases, the number of samples used for the outer average also increases, resulting in higher accuracy for the outer expectation, whereas the number of samples for the inner average decreases, resulting in lower accuracy for the inner expectation. The previous method by the authors \cite{hironaka2023efficient}, roughly speaking, picks up only one summand with $d=m/2$ without the second term of \eqref{eq:proposed_estimator} to estimate the nested expectation, i.e., $\hat{P}^{(m/2,m/2)}$, which aims to balance the accuracies for approximating the outer and inner expectations. Instead, in this paper, we use the whole hierarchy of the approximation levels and correct the estimator based on the idea of sparse grid quadrature as follows.

To understand how our estimator is connected to sparse grid quadrature, let us consider the second term of \eqref{eq:proposed_estimator}, which is a sum over $d$ from 1 to $m$. We denote each summand by
\begin{align*}
    \hat{Q}^{(d,m-d)} & := \frac{1}{2^d}\sum_{\mathbf{u} \in \{0, 1\}^d}f\left(\frac{1}{2^{m-d}}\sum_{i \in B'^d_{\mathbf{u}}} X^{(i)}\right)\\
    & \: = \frac{1}{2^{d-1}}\sum_{\mathbf{u} \in \{0, 1\}^{d-1}}\frac{1}{2}\sum_{p\in \{0,1\}}f\left(\frac{1}{2^{m-d}}\sum_{i \in B'^d_{\mathbf{u}{:}p}} X^{(i)}\right).
\end{align*}
Although $\hat{Q}^{(d,m-d)}$ can be computed similarly to $\hat{P}^{(d,m-d)}$, it is important to notice that, for each $\mathbf{u}\in \{0,1\}^{d-1}$, $B'^d_{\mathbf{u}{:}0}$ and $B'^d_{\mathbf{u}{:}1}$ partition the samples in $B_{\mathbf{u}}^{d-1}$ depending only on the index $i$ of the samples $\{(X^{(i)}, Y^{(i)})\}_{1 \leq i \leq 2^m}$. Due to the assumption that $\{(X^{(i)}, Y^{(i)})\}_{1 \leq i \leq 2^m}$ is a set of i.i.d.\ random samples from the joint distribution, the statistical property for the averages over $B'^d_{\mathbf{u}{:}0}$ and $B'^d_{\mathbf{u}{:}1}$ remains the same. Therefore, $\hat{Q}^{(d,m-d)}$ plays a similar role to $\hat{P}^{(d-1,m-d)}$ in approximating the nested expectation. Following the idea of sparse grid quadrature \cite{smolyak1963quadrature,Gerstner1998Numerical,bungartz2004sparse}, we aim to approximate the nested expectation, with the accuracy which $\hat{P}^{(m,m)}$ can achieve, by
\[ \hat{P}^{(0,m)}+\sum_{d=1}^{m}\left(\hat{P}^{(d,m-d)}-\hat{Q}^{(d,m-d)}\right),\]
which is nothing but our proposed estimator in \eqref{eq:proposed_estimator}.

\begin{remark}\label{rem:dom}
    In Algorithm~\ref{alg:proposed}, constructing $B^{d}_{\mathbf{u}}$ for $d\in \{0,\dots,m\}$ and $\mathbf{u}\in \{0,1\}^d$ requires recursive sorting of the samples based on the value of $Y^{(i)}$. For each $d$, the necessary computational cost for this sorting process is of $O(N\log N)$, where $N=2^m$ denotes the total number of samples. Therefore, the total cost will be of $O(mN\log N)=O(N(\log N)^2)$, which is almost linear in $N$.

    In the rest of this paper, we assume that the cost of recursive sorting is negligible, and the dominant computational parts in Algorithm~\ref{alg:proposed} are generating random samples from the joint distribution $\rho(X,Y)$ and evaluating the function $f$, both of which cost linearly with respect to the total number of samples used. Therefore, we adopt the total number of samples used in estimation as an objective measure of computational cost when comparing the performance of different estimation methods.
\end{remark}

%% file: 30-theoretical_th.tex
\section{Theoretical Analysis}\label{sec:theoretical}

In this section, we provide an upper bound on the MSE of our proposed method (Algorithm~\ref{alg:proposed}). We begin by introducing our notation for variance. Throughout this section, given a random vector $Z$, we define the variance of $Z$ as the trace of the covariance matrix of $Z$ and denote it by $\mathbb{V}_{\rho(Z)}Z$. This notation also applies to the case of conditional variance, where we use an appropriate subscript to represent the underlying conditional distribution. Moreover, whenever the underlying probability distribution is clear from the context, we omit the subscript from $\mathbb{E}$ or $\mathbb{V}$ or simply denote the underlying random variables for the subscript instead of the probability distribution.

As the main theoretical result of this paper, we show the following theorem.
\begin{theorem}\label{thm:main}
    Let $m\in \mathbb{N}$ and $N = 2^m$.
   Assume that $\sup_{Y}\mathbb{V}_{X|Y} X<\infty$, and there exist positive constants $\alpha_1, \alpha_2$ such that
    \begin{enumerate}
        \item $\left| f\left( x \right) - f\left( x + y \right) \right| \leq \alpha_1 \left\| y \right\|_2$, and
        \item $\left|f\left( x \right)-f\left( x + y \right)-f\left( x + y' \right)+f\left( x + y + y' \right)\right|\leq \alpha_2 \left\| y \right\|_2 \left\| y' \right\|_2$,
    \end{enumerate}
    for all $x, y, y' \in \mathbb{R}^J$.
    Moreover, assume that there exist a positive constant $\beta$ and functions $\mathcal{T}_k: \mathbb{R}\to [0,1]$ for all $k=1,\ldots,K$ such that, for any pair of $Y_1$ and $Y_2$, there is a probability distribution $\rho_{Y_1, Y_2}$ which satisfies, for all $(X_1, X_2) \sim \rho_{Y_1, Y_2}$,
    \begin{enumerate}\setcounter{enumi}{2}
        \item $X_1 \sim \rho(X\mid Y_1)$, $X_2 \sim \rho(X\mid Y_2)$, and
        \item $\left\| X_1 - X_2 \right\|_2 \leq \beta \left\| \mathcal{T}(Y_1) - \mathcal{T}(Y_2) \right\|_2$,
    \end{enumerate}
    where we write $\mathcal{T}(Y)=(\mathcal{T}_1(Y_1),\dots,\mathcal{T}_K(Y_K)).$
    Then we have
    \begin{align}
        \mathbb{E}(I - \hat{I})^2
    =
        O\left( \frac{(\log N)^2}{N^{1/K}} \right).
    \end{align}
\end{theorem}

Before going into the proof of the theorem, let us make some remarks.
\begin{remark}\label{rem:T_k}
The assumption of the existence of functions $\mathcal{T}_k$ may appear artificial. To argue that this is not the case, let us consider the following change of variables: for $k=1,\dots,K$, let $F_k$ denote the cumulative marginal distribution function of the outer variable $Y_k$ and let $Z_k=F_k(Y_k)$. This transformation ensures that $0 \leq Z_k \leq 1$. Therefore, a natural choice for $\mathcal{T}_k$ is $\mathcal{T}_k=F_k$.

Furthermore, note that the third item of Theorem~\ref{thm:main} implies that the random variables $X_1$ and $X_2$ following the joint distribution $\rho_{Y_1, Y_2}$ marginally follow the conditional distributions $\rho(X\mid Y_1)$ and $\rho(X\mid Y_2)$, respectively. The fourth item of Theorem~\ref{thm:main} assumes that the inner conditional distribution can be coupled well depending on the closeness of $Y_1$ and $Y_2$. It is important to note that we do not require such coupling in our algorithm; we just need to assume this property in the proof. These assumptions lead to the following bound:
\[ \left\|\mathbb{E}_{X|Y}X - \mathbb{E}_{X'|Y'}X'\right\|_2 \leq \beta\left\| \mathcal{T}(Y) - \mathcal{T}(Y') \right\|_2,\]
which holds for all $Y, Y' \in \mathbb{R}^K$.
\end{remark}

\begin{remark}
If the function $f$ satisfies the first condition, it also satisfies
\begin{align*}
    \left|f\left( x \right)-f\left( x + y \right)-f\left( x + y' \right)+f\left( x + y + y' \right)\right| \leq 2\alpha_1\min\left\{ \left\| y \right\|_2,\left\| y' \right\|_2\right\}\leq 2\alpha_1\left( \left\| y \right\|_2 \left\| y' \right\|_2\right)^{1/2}.
\end{align*}
Therefore, it is clear that the second condition imposes a stronger property for $f$.

Theorem~\ref{thm:main} does not necessarily cover the problem setting for estimating EVSI. In fact, for the case where $f=\max_{d\in D}$ for a set $D$ with finite cardinality $J$, it can be checked that the first condition on $f$ is satisfied:
\begin{align*}
    \left| \max_{1\leq j\leq J}x_j-\max_{1\leq j\leq J}(x_j+y_j)\right|\leq \max_{1\leq j\leq J}\left| x_j-(x_j+y_j)\right|\leq \sum_{j=1}^{J}\left| y_j\right|\leq \sqrt{J}\, \|y\|_2,
\end{align*}
where the last inequality follows from Cauchy–Schwarz inequality. However, the second condition on $f$ is not satisfied. Therefore, discussing whether the conditions given in Theorem~\ref{thm:main} can be weakened such that the theorem covers the EVSI setting remains an open problem.
\end{remark}

As a preparation for proving Theorem~\ref{thm:main}, we present the following lemma.
\begin{lemma}\label{lem}
    Let $m \in \mathbb{N}$, and let $\{(X^{(i)}, Y^{(i)})\}_{1 \leq i \leq 2^m}$ be a set of i.i.d.\ random samples from the joint distribution $\rho(X,Y)$. For $d \in \{0, \dots, m\}$ and $\mathbf{u} \in \{0,1\}^d$, let $B^{d}_{\mathbf{u}}$ be constructed as described in Algorithm~\ref{alg:proposed}. Then, assuming the existence of functions $\mathcal{T}_k: \mathbb{R}\to [0,1]$ for all $k=1,\dots, K$  we have for any $d \in \{0, \dots, m\}$
    \begin{align}
        \frac{1}{2^d}
        \sum_{\mathbf{u} \in \{0, 1\}^d}
        \max_{i, j \in B^d_{\mathbf{u}}}
        \left\|
            \mathcal{T}(Y^{(i)})
        -
            \mathcal{T}(Y^{(j)})
        \right\|_2^2
    \leq
        \frac{2K}{2^{d/K}}.
    \end{align}
\end{lemma}

\begin{proof}
    Let $N=2^m$. For $d \in \{0, \dots, m\}$ and $k\in \{1,\dots, K\}$, define
    \begin{align}
        W_{d, k}
    =
        \frac{1}{2^d}
        \sum_{\mathbf{u} \in \{0, 1\}^d}
        \left(
            \max_{i \in B^d_{\mathbf{u}}}
            \mathcal{T}_k(Y_k^{(i)})
        -
            \min_{i \in B^d_{\mathbf{u}}}
            \mathcal{T}_k(Y_k^{(i)})
        \right).
    \end{align}

    If $d = 0$, it is easy to see from the definition $B^0_{\emptyset}:=\{1,\dots,2^m\}$ that
    \begin{align}
        W_{0, k}
    =
        \max_{i \in B^0_{\emptyset}}
        \mathcal{T}_k(Y_k^{(i)})
    -
        \min_{i \in B^0_{\emptyset}}
        \mathcal{T}_k(Y_k^{(i)})
    =
        \max_{ i=1, \dots, 2^m }
        \mathcal{T}_k(Y_k^{(i)})
    -
        \min_{ i=1, \dots, 2^m }
        \mathcal{T}_k(Y_k^{(i)})
    \leq
        1
    \end{align}
    holds for all $k\in\{1,\dots,K\}$.

    Otherwise if $1 \leq d \leq m$, it follows from the construction of $B^{d}_{\mathbf{u}}$'s that, for $k = (d-1 \bmod K) +1 \in \{1,\dots,K\}$, we have
    \begin{align}
    &
        \min_{i \in B_{\mathbf{u}}^{d-1}}
        \mathcal{T}_k(Y_k^{(i)})
    =
        \min_{i \in B_{\mathbf{u}{:}0}^d}
        \mathcal{T}_k(Y_k^{(i)})
    \leq
        \max_{i \in B_{\mathbf{u}{:}0}^d}
        \mathcal{T}_k(Y_k^{(i)})
    \\ &\quad \quad \leq
        \min_{i \in B_{\mathbf{u}{:}1}^d}
        \mathcal{T}_k(Y_k^{(i)})
    \leq
        \max_{i \in B_{\mathbf{u}{:}1}^d}
        \mathcal{T}_k(Y_k^{(i)})
    =
        \max_{i \in B_{\mathbf{u}}^{d-1}}
        \mathcal{T}_k(Y_k^{(i)})
    \end{align}
    for all $\mathbf{u} \in \{0, 1\}^{d-1}$. Therefore, we obtain
    \begin{align}
        W_{d, k} &=
        \frac{1}{2^d}
        \sum_{\mathbf{u} \in \{0, 1\}^d}
        \left(
            \max_{i \in B^d_{\mathbf{u}}}
            \mathcal{T}_k(Y_k^{(i)})
        -
            \min_{i \in B^d_{\mathbf{u}}}
            \mathcal{T}_k(Y_k^{(i)})
        \right)
    \\ &=
        \frac{1}{2^d}
        \sum_{\mathbf{u} \in \{0, 1\}^{d-1}}
        \left(
            \left(
                \max_{i \in B^d_{\mathbf{u}{:}0}}
                \mathcal{T}_k(Y_k^{(i)})
            -
                \min_{i \in B^d_{\mathbf{u}{:}0}}
                \mathcal{T}_k(Y_k^{(i)})
            \right)
        +
            \left(
                \max_{i \in B^d_{\mathbf{u}{:}1}}
                \mathcal{T}_k(Y_k^{(i)})
            -
                \min_{i \in B^d_{\mathbf{u}{:}1}}
                \mathcal{T}_k(Y_k^{(i)})
            \right)
        \right)
    \\ &\leq
        \frac{1}{2^d}
        \sum_{\mathbf{u} \in \{0, 1\}^{d-1}}
        \left(
            \max_{i \in B^{d-1}_{\mathbf{u}}}
            \mathcal{T}_k(Y_k^{(i)})
        -
            \min_{i \in B^{d-1}_{\mathbf{u}}}
            \mathcal{T}_k(Y_k^{(i)})
        \right) =
        \frac{1}{2}
        W_{d-1, k}.
    \end{align}
    For $k \neq (d-1 \bmod K) +1 $, we have
    \begin{align}
    \quad
        \min_{i \in B_{\mathbf{u}}^{d-1}}
        \mathcal{T}_k(Y_k^{(i)})
    \leq
        \min_{i \in B_{\mathbf{u}{:}p}^d}
        \mathcal{T}_k(Y_k^{(i)})
    \leq
        \max_{i \in B_{\mathbf{u}{:}p}^d}
        \mathcal{T}_k(Y_k^{(i)})
    \leq
        \max_{i \in B_{\mathbf{u}}^{d-1}}
        \mathcal{T}_k(Y_k^{(i)})
    \end{align}
    for all $\mathbf{u} \in \{0, 1\}^{d-1}$ and $p \in \{0, 1\}$. Therefore, we obtain
    \begin{align}
     W_{d, k} &=
        \frac{1}{2^d}
        \sum_{\mathbf{u} \in \{0, 1\}^d}
        \left(
            \max_{i \in B^d_{\mathbf{u}}}
            \mathcal{T}_k(Y_k^{(i)})
        -
            \min_{i \in B^d_{\mathbf{u}}}
            \mathcal{T}_k(Y_k^{(i)})
        \right)
    \\ &=
        \frac{1}{2^d}
        \sum_{\mathbf{u} \in \{0, 1\}^{d-1}}
        \left(
            \left(
                \max_{i \in B^d_{\mathbf{u}{:}0}}
                \mathcal{T}_k(Y_k^{(i)})
            -
                \min_{i \in B^d_{\mathbf{u}{:}0}}
                \mathcal{T}_k(Y_k^{(i)})
            \right)
        +
            \left(
                \max_{i \in B^d_{\mathbf{u}{:}1}}
                \mathcal{T}_k(Y_k^{(i)})
            -
                \min_{i \in B^d_{\mathbf{u}{:}1}}
                \mathcal{T}_k(Y_k^{(i)})
            \right)
        \right)
    \\ &\leq
        \frac{2}{2^d}
        \sum_{\mathbf{u} \in \{0, 1\}^{d-1}}
        \left(
            \max_{i \in B^{d-1}_{\mathbf{u}}}
            \mathcal{T}_k(Y_k^{(i)})
        -
            \min_{i \in B^{d-1}_{\mathbf{u}}}
            \mathcal{T}_k(Y_k^{(i)})
        \right) =
        W_{d-1, k}.
    \end{align}
    Using these results, it can be concluded that, for any $k\in \{1,\dots,K\}$, we have
    \begin{align}
        W_{d, k} \leq \begin{cases} 1 & \text{if $d < K$,}\\ W_{d-K, k}/2 & \text{otherwise.} \end{cases}
    \end{align}
    so that it holds that
    \begin{align}
            W_{d, k}
        \leq
            \frac{1}{2^{\lfloor d/K \rfloor}}
            W_{d-K\lfloor d/K \rfloor, k}
        \leq
            \frac{1}{2^{\lfloor d/K \rfloor}}
        \leq
            \frac{2}{2^{d/K}},
    \end{align}
    for any $1\leq d\leq m$. Note that this upper bound on $W_{d, k}$ also applies to the case $d=0$.

    Recalling the assumption that $0 \leq \mathcal{T}_k(Y_k) \leq 1$ for any $Y_k$ for all $k=1,\dots, K$, we have for any $d \in \{0, \dots, m\}$
    \begin{align}
        \frac{1}{2^d}
        \sum_{\mathbf{u} \in \{0, 1\}^d}
        \max_{i, j \in B^d_{\mathbf{u}}}
        \left\|
            \mathcal{T}(Y^{(i)})
        -
            \mathcal{T}(Y^{(j)})
        \right\|_2^2 &\leq
        \frac{1}{2^d}
        \sum_{\mathbf{u} \in \{0, 1\}^d}
        \max_{i, j \in B^d_{\mathbf{u}}}
        \left\|
            \mathcal{T}(Y^{(i)})
        -
            \mathcal{T}(Y^{(j)})
        \right\|_1
    \\ &\leq
        \frac{1}{2^d}
        \sum_{k=1}^{K}
        \sum_{\mathbf{u} \in \{0, 1\}^d}
        \max_{i, j \in B^d_{\mathbf{u}}}
        \left|
            \mathcal{T}_k(Y_k^{(i)})
        -
            \mathcal{T}_k(Y_k^{(j)})
        \right|
    \\ &=
        \frac{1}{2^d}
        \sum_{k=1}^{K}
        \sum_{\mathbf{u} \in \{0, 1\}^d}
        \left(
            \max_{i \in B^d_{\mathbf{u}}}
            \mathcal{T}_k(Y_k^{(i)})
        -
            \min_{i \in B^d_{\mathbf{u}}}
            \mathcal{T}_k(Y_k^{(i)})
        \right)
    \\ &=
        \sum_{k=1}^{K}
        W_{d, K} \leq
        \frac{2K}{2^{d/K}},
    \end{align}
    which completes the proof of the lemma.
\end{proof}

Let us now prove Theorem~\ref{thm:main}.

\begin{proof}[Proof of Theorem~\ref{thm:main}]
    Together with the notation used in Section~\ref{sec:method}, we further introduce the following notation in this proof:
    \begin{align*}
        \hat{P}_{\mathrm{in}}^{(d,m-d)} & = \frac{1}{2^d}\sum_{\mathbf{u} \in \{0, 1\}^d}f\left(\frac{1}{2^{m-d}}\sum_{i \in B^d_{\mathbf{u}}}\mathbb{E}_{X|Y^{(i)}}X\right), \\
        \hat{Q}_{\mathrm{in}}^{(d,m-d)} & = \frac{1}{2^d}\sum_{\mathbf{u} \in \{0, 1\}^d}f\left(\frac{1}{2^{m-d}}\sum_{i \in B'^d_{\mathbf{u}}}\mathbb{E}_{X|Y^{(i)}}X\right), \\
        \hat{P}_{\mathrm{out}}^{(d,m-d)} & = \frac{1}{2^d}\sum_{\mathbf{u} \in \{0, 1\}^d}\mathbb{E}_{\{X^{(i)}|Y^{(i)}\}}f\left(\frac{1}{2^{m-d}}\sum_{i \in B^d_{\mathbf{u}}}X^{(i)}\right), \\
        \hat{Q}_{\mathrm{out}}^{(d,m-d)} & = \frac{1}{2^d}\sum_{\mathbf{u} \in \{0, 1\}^d}\mathbb{E}_{\{X^{(i)}|Y^{(i)}\}}f\left(\frac{1}{2^{m-d}}\sum_{i \in B'^d_{\mathbf{u}}}X^{(i)}\right).
    \end{align*}
    Then, by using a trivial inequality $(a+b+c)^2\leq 3(a^2+b^2+c^2)$ for any $a,b,c\in \mathbb{R}$, which follows immediately from Cauchy–Schwarz inequality, we have
    \begin{align}
    &\quad
        \mathbb{E}(I - \hat{I})^2
    \\ &=
        \mathbb{E}
        \left(
            \mathbb{E}_Y
            f\left(
                \mathbb{E}_{X|Y}X
            \right)
        -
            \left(
                \sum_{d=0}^{m}
                \hat{P}^{(d,m-d)}
            -
                \sum_{d=1}^{m}
                \hat{Q}^{(d,m-d)}
            \right)
        \right)^2
    \\ &\leq
        3
        \mathbb{E}_{\{Y^{(i)}\}}
        \left(
            \mathbb{E}_Y
            f\left(
                \mathbb{E}_{X|Y}X
            \right)
        -
            \left(
                \sum_{d=0}^{m}
                \hat{P}_{\mathrm{in}}^{(d,m-d)}
            -
                \sum_{d=1}^{m}
                \hat{Q}_{\mathrm{in}}^{(d,m-d)}
            \right)
        \right)^2
    \\ &\quad\quad
    +
        3
        \mathbb{E}_{\{Y^{(i)}\}}
        \left(
            \left(
                \sum_{d=0}^{m}
                \hat{P}_{\mathrm{in}}^{(d,m-d)}
            -
                \sum_{d=1}^{m}
                \hat{Q}_{\mathrm{in}}^{(d,m-d)}
            \right)
        -
            \left(
                \sum_{d=0}^{m}
                \hat{P}_{\mathrm{out}}^{(d,m-d)}
            -
                \sum_{d=1}^{m}
                \hat{Q}_{\mathrm{out}}^{(d,m-d)}
            \right)
        \right)^2
    \\ &\quad\quad
    +
        3
        \mathbb{E}_{\{(X^{(i)},Y^{(i)})\}}
        \left(
            \left(
                \sum_{d=0}^{m}
                \hat{P}_{\mathrm{out}}^{(d,m-d)}
            -
                \sum_{d=1}^{m}
                \hat{Q}_{\mathrm{out}}^{(d,m-d)}
            \right)
        -
            \left(
                \sum_{d=0}^{m}
                \hat{P}^{(d,m-d)}
            -
                \sum_{d=1}^{m}
                \hat{Q}^{(d,m-d)}
            \right)
        \right)^2.\label{eq:main_theorem_three_terms}
    \end{align}
In the following, we derive an upper bound on each of the three terms of \eqref{eq:main_theorem_three_terms}.

\underline{A bound on the first term of \eqref{eq:main_theorem_three_terms}:}
We start by considering the first term. Applying Cauchy–Schwarz inequality for sums, we obtain
    \begin{align}
    &\quad
        \mathbb{E}_{\{Y^{(i)}\}}
        \left(
            \mathbb{E}_Y
            f\left(
                \mathbb{E}_{X|Y}X
            \right)
        -
            \left(
                \sum_{d=0}^{m}
                \hat{P}_{\mathrm{in}}^{(d,m-d)}
            -
                \sum_{d=1}^{m}
                \hat{Q}_{\mathrm{in}}^{(d,m-d)}
            \right)
        \right)^2
    \\ &=\mathbb{E}_{\{Y^{(i)}\}}
        \left(
            \mathbb{E}_Y
            f\left(
                \mathbb{E}_{X|Y}X
            \right)
        -
            \hat{P}_{\mathrm{in}}^{(m,0)}
        -
            \sum_{d=1}^{m}
                \left(
                \hat{P}_{\mathrm{in}}^{(d-1,m-d+1)}
            -
                \hat{Q}_{\mathrm{in}}^{(d,m-d)}
            \right)
        \right)^2
    \\ &\leq
        (m+1)
        \mathbb{E}_{\{Y^{(i)}\}}
        \left(
            \mathbb{E}_Y
            f\left(
                \mathbb{E}_{X|Y} X
            \right)
        -
            \hat{P}_{\mathrm{in}}^{(m,0)}
        \right)^2
    +
        (m+1)
        \sum_{d=1}^{m}
        \mathbb{E}_{\{Y^{(i)}\}}
        \left(
                \hat{P}_{\mathrm{in}}^{(d-1,m-d+1)}
            -
                \hat{Q}_{\mathrm{in}}^{(d,m-d)}
        \right)^2.
    \end{align}
For the first term on the right-most side above, we have
    \begin{align}
        \mathbb{E}_{\{Y^{(i)}\}}
        \left(
            \mathbb{E}_Y
            f\left(
                \mathbb{E}_{X|Y} X
            \right)
        -
             \hat{P}_{\mathrm{in}}^{(m,0)}
        \right)^2
        =
        \mathbb{E}_{\{Y^{(i)}\}}
        \left(
            \mathbb{E}_Y
            f\left(
                \mathbb{E}_{X|Y} X
            \right)
        -
            \frac{1}{2^m}
            \sum_{i=1}^{2^m}
            f\left(
                \mathbb{E}_{X|Y^{(i)}} X
            \right)
        \right)^2
        =
        \frac{1}{2^m}
        \mathbb{V}_{Y}
        f\left(
            \mathbb{E}_{X|Y} X
        \right).
    \end{align}
Regarding each summand of the second term, applying Jensen's inequality and using the first condition on $f$, we have
    \begin{align}
    &\quad
        \mathbb{E}_{\{Y^{(i)}\}}
        \left(
                \hat{P}_{\mathrm{in}}^{(d-1,m-d+1)}
            -
                \hat{Q}_{\mathrm{in}}^{(d,m-d)}
        \right)^2
    \\ &=
        \mathbb{E}_{\{Y^{(i)}\}}
        \left(
        \frac{1}{2^{d-1}}
        \sum_{\mathbf{u} \in \{0, 1\}^{d-1}}
            \left(
            f\left(
                \frac{1}{2^{m-d+1}}
                \sum_{i \in B^{d-1}_{\mathbf{u}}}
                \mathbb{E}_{X|Y^{(i)}} X
            \right)
        -
            \frac{1}{2}
            \sum_{p \in \{0, 1\}}
            f\left(
                \frac{1}{2^{m-d}}
                \sum_{i \in B'^{d}_{\mathbf{u}{:}p}}
                \mathbb{E}_{X|Y^{(i)}} X
            \right)
            \right)
        \right)^2
    \\ & \leq
        \mathbb{E}_{\{Y^{(i)}\}}
        \frac{1}{2^{d}}
        \sum_{\mathbf{u} \in \{0, 1\}^{d-1}}\sum_{p \in \{0, 1\}}
        \left(
            f\left(
                \frac{1}{2^{m-d+1}}
                \sum_{i \in B^{d-1}_{\mathbf{u}}}
                \mathbb{E}_{X|Y^{(i)}} X
            \right)
        -
            f\left(
                \frac{1}{2^{m-d}}
                \sum_{i \in B'^{d}_{\mathbf{u}{:}p}}
                \mathbb{E}_{X|Y^{(i)}} X
            \right)
        \right)^2
    \\ & \leq
        \alpha_1^2\mathbb{E}_{\{Y^{(i)}\}}
        \frac{1}{2^{d}}
        \sum_{\mathbf{u} \in \{0, 1\}^{d-1}}\sum_{p \in \{0, 1\}}
        \left\|
                \frac{1}{2^{m-d+1}}
                \sum_{i \in B^{d-1}_{\mathbf{u}}}
                \mathbb{E}_{X|Y^{(i)}} X
        -
                \frac{1}{2^{m-d}}
                \sum_{i \in B'^{d}_{\mathbf{u}{:}p}}
                \mathbb{E}_{X|Y^{(i)}} X
        \right\|_2^2.
    \end{align}
    By construction of $B'^{d}_{\mathbf{u}{:}p}$ and because of the exchangeability of the samples $\{Y^{(i)}\}_i$, the set $\{Y^{(i)}\}_{i \in B'^{d}_{\mathbf{u}{:}p}}$ can be regarded as a set of the samples randomly resampled from $\{Y^{(i)}\}_{i \in B^{d-1}_{\mathbf{u}}}$ without replacement. This implies that each summand $\mathbb{E}_{X|Y^{(i)}} X$ for $i \in B'^{d}_{\mathbf{u}{:}p}$ is an unbiased estimator of its average over $B^{d-1}_{\mathbf{u}}$. Therefore, the expectation of each summand above corresponds to the variance of an unbiased estimator that uses $2^{m-d}$ samples (without replacement) to estimate the expectation of a finite population of size $2^{m-d+1}$. It follows from this argument that
    \begin{align}
    &\quad
        \mathbb{E}_{\{Y^{(i)}\}}
        \left(
                \hat{P}_{\mathrm{in}}^{(d-1,m-d+1)}
            -
                \hat{Q}_{\mathrm{in}}^{(d,m-d)}
        \right)^2
    \\ & \leq
        \alpha_1^2\mathbb{E}_{\{Y^{(i)}\}}
        \frac{1}{2^{d}}\sum_{\mathbf{u} \in \{0, 1\}^{d-1}}\sum_{p \in \{0, 1\}}\frac{2^{m-d+1}-2^{m-d}}{2^{m-d+1}-1}\cdot
        \frac{1}{2^{m-d}}
    \\ &\quad\quad \times
        \frac{1}{2^{m-d+1}}
            \sum_{i \in B^{d-1}_{\mathbf{u}}}
            \left\|
                \frac{1}{2^{m-d+1}}
                \sum_{j \in B^{d-1}_{\mathbf{u}}}
                \mathbb{E}_{X|Y^{(j)}} X
            -
                \mathbb{E}_{X|Y^{(i)}} X
            \right\|_2^2
    \\ &=
        \alpha_1^2
        \mathbb{E}_{\{Y^{(i)}\}}
        \frac{1}{2^m-2^{d-1}}
        \sum_{\mathbf{u} \in \{0, 1\}^{d-1}}
            \frac{1}{2^{m-d+1}}
            \sum_{i \in B^{d-1}_{\mathbf{u}}}
            \left\|
                \frac{1}{2^{m-d+1}}
                \sum_{j \in B^{d-1}_{\mathbf{u}}}
                \mathbb{E}_{X|Y^{(j)}} X
            -
                \mathbb{E}_{X|Y^{(i)}} X
            \right\|_2^2
    \\ &\leq
        \alpha_1^2
        \mathbb{E}_{\{Y^{(i)}\}}
        \frac{1}{2^{m-1}}
        \sum_{\mathbf{u} \in \{0, 1\}^{d-1}}
        \max_{i, j \in B^{d-1}_{\mathbf{u}}}
        \left\|
            \mathbb{E}_{X|Y^{(i)}} X
        -
            \mathbb{E}_{X|Y^{(j)}} X
        \right\|_2^2
    \\ &\leq
        \alpha_1^2\beta^2
        \mathbb{E}_{\{Y^{(i)}\}}
        \frac{1}{2^{m-1}}
        \sum_{\mathbf{u} \in \{0, 1\}^{d-1}}
        \max_{i, j \in B^{d-1}_{\mathbf{u}}}
        \left\|
            \mathcal{T}(Y^{(i)})
        -
            \mathcal{T}(Y^{(j)})
        \right\|_2^2
    \\ &\leq
        \alpha_1^2\beta^2
        \frac{2^{d-1}}{2^{m-1}}\cdot
        \frac{2K}{2^{(d-1)/K}}
        \leq
        \alpha_1^2
        \beta^2
        \frac{2K}{2^{(m-1)/K}},
    \end{align}
    in which we have used the assumption on the coupling of $X$ (see also Remark~\ref{rem:T_k}) and Lemma~\ref{lem}. In this way, the first term of \eqref{eq:main_theorem_three_terms} is bounded above by
    \begin{align}
    &\quad
        \mathbb{E}_{\{Y^{(i)}\}}
        \left(
            \mathbb{E}_Y
            f\left(
                \mathbb{E}_{X|Y} X
            \right)
        -
            \left(
                \sum_{d=0}^{m}
                \hat{P}_{\mathrm{in}}^{(d,m-d)}
            -
                \sum_{d=1}^{m}
                \hat{Q}_{\mathrm{in}}^{(d,m-d)}
            \right)
        \right)^2
    \\ &\leq
        (m+1)
        \mathbb{E}_{\{Y^{(i)}\}}
        \left(
            \mathbb{E}_Y
            f\left(
                \mathbb{E}_{X|Y} X
            \right)
        -
            \hat{P}_{\mathrm{in}}^{(m,0)}
        \right)^2
    +
        (m+1)
        \sum_{d=1}^{m}
        \mathbb{E}_{\{Y^{(i)}\}}
        \left(
                \hat{P}_{\mathrm{in}}^{(d-1,m-d+1)}
            -
                \hat{Q}_{\mathrm{in}}^{(d,m-d)}
        \right)^2
    \\ &\leq
        \frac{m+1}{2^m}
        \mathbb{V}_{Y}
        f\left(
            \mathbb{E}_{X|Y} X
        \right)
    +
        (m+1) \sum_{d=1}^{m}\alpha_1^2
        \beta^2
        \frac{2K}{2^{(m-1)/K}}
    =
        O\left(
            \frac{m^2}{2^{m/K}}
        \right).
    \end{align}

\underline{A bound on the second term of \eqref{eq:main_theorem_three_terms}:}
    Next we show a bound on the second term of \eqref{eq:main_theorem_three_terms}. Using Cauchy–Schwarz inequality for sums, we have
    \begin{align}
    &\quad
        \mathbb{E}_{\{Y^{(i)}\}}
        \left(
            \left(
                \sum_{d=0}^{m}
                \hat{P}_{\mathrm{in}}^{(d,m-d)}
            -
                \sum_{d=1}^{m}
                \hat{Q}_{\mathrm{in}}^{(d,m-d)}
            \right)
        -
            \left(
                \sum_{d=0}^{m}
                \hat{P}_{\mathrm{out}}^{(d,m-d)}
            -
                \sum_{d=1}^{m}
                \hat{Q}_{\mathrm{out}}^{(d,m-d)}
            \right)
        \right)^2
    \\ &=
        \mathbb{E}_{\{Y^{(i)}\}}
        \left(
            \hat{P}_{\mathrm{in}}^{(0,m)}-\hat{P}_{\mathrm{out}}^{(0,m)}
        +
            \sum_{d=1}^{m}
               \left( \hat{P}_{\mathrm{in}}^{(d,m-d)}
            -
                \hat{P}_{\mathrm{out}}^{(d,m-d)}
                \right)
        -
            \sum_{d=1}^{m}
               \left( \hat{Q}_{\mathrm{in}}^{(d,m-d)}
            -
                \hat{Q}_{\mathrm{out}}^{(d,m-d)}
                \right)
        \right)^2
    \\ &\leq
       (m+1) \mathbb{E}_{\{Y^{(i)}\}}
        \left(
            \hat{P}_{\mathrm{in}}^{(0,m)}-\hat{P}_{\mathrm{out}}^{(0,m)}
        \right)^2
    \\ &\quad
        +
       (m+1) \sum_{d=1}^{m}\mathbb{E}_{\{Y^{(i)}\}}
        \left(
                \hat{P}_{\mathrm{in}}^{(d,m-d)}
            -
                \hat{P}_{\mathrm{out}}^{(d,m-d)}
        -
                \hat{Q}_{\mathrm{in}}^{(d,m-d)}
            +
                \hat{Q}_{\mathrm{out}}^{(d,m-d)}
        \right)^2.
    \end{align}
    The first term on the right-most side above can be bounded above by
    \begin{align}
        \mathbb{E}_{\{Y^{(i)}\}}
        \left(
            \hat{P}_{\mathrm{in}}^{(0,m)}-\hat{P}_{\mathrm{out}}^{(0,m)}
        \right)^2
        &=
        \mathbb{E}_{\{Y^{(i)}\}}
        \left(
            f\left(
                \frac{1}{2^m}
                \sum_{i=1}^{2^m}
                \mathbb{E}_{X|Y^{(i)}} X
            \right)
        -
            \mathbb{E}_{\{X^{(i)}|Y^{(i)}\}}
            f\left(
                \frac{1}{2^m}
                \sum_{i=1}^{2^m} X^{(i)}
            \right)
        \right)^2
    \\ &\leq
        \mathbb{E}_{\{(X^{(i)},Y^{(i)})\}}
        \left(
            f\left(
                \frac{1}{2^m}
                \sum_{i=1}^{2^m}
                \mathbb{E}_{X|Y^{(i)}} X
            \right)
        -
            f\left(
                \frac{1}{2^m}
                \sum_{i=1}^{2^m} X^{(i)}
            \right)
        \right)^2
    \\ &\leq
        \alpha_1^2
        \mathbb{E}_{\{(X^{(i)},Y^{(i)})\}}
        \left\|
            \frac{1}{2^m}
            \sum_{i=1}^{2^m}
            \mathbb{E}_{X|Y^{(i)}} X
        -
            \frac{1}{2^m}
            \sum_{i=1}^{2^m} X^{(i)}
        \right\|_2^2
    \\ &=
        \alpha_1^2
        \mathbb{E}_{\{Y^{(i)}\}}
        \frac{1}{2^{2m}}
        \sum_{i=1}^{2^m}
        \mathbb{V}_{X|Y^{(i)}} X
        =
        \frac{\alpha_1^2}{2^m}
        \mathbb{E}_{Y}
        \mathbb{V}_{X|Y} X
    \end{align}
    where the first inequality follows from Jensen's inequality and the second inequality follows from the first condition of $f$. Regarding each summand of the second term, let us introduce a bijection $\sigma: \{1,\dots,2^m\} \to \{1,\dots,2^m\}$, which is also locally bijective between $B^{d}_{\mathbf{u}}$ and $B'^{d}_{\mathbf{u}}$, i.e., $\sigma(B^{d}_{\mathbf{u}})=B'^{d}_{\mathbf{u}}$ holds for all $\mathbf{u}\in \{0,1\}^d$. Then, by using the first and second conditions on $f$, the third item of the assumptions in the theorem (see also Remark~\ref{rem:T_k}) and Cauchy–Schwarz inequality, we have
    \begin{align}
    &\quad
        \mathbb{E}_{\{Y^{(i)}\}}
        \left(
                \hat{P}_{\mathrm{in}}^{(d,m-d)}
            -
                \hat{P}_{\mathrm{out}}^{(d,m-d)}
        -
                \hat{Q}_{\mathrm{in}}^{(d,m-d)}
            +
                \hat{Q}_{\mathrm{out}}^{(d,m-d)}
        \right)^2
    \\ &=
        \mathbb{E}_{\{Y^{(i)}\}}
        \left(
            \frac{1}{2^d}
            \sum_{\mathbf{u} \in \{0, 1\}^d}
            \left(
                f\left(
                    \frac{1}{2^{m-d}}
                    \sum_{i \in B^d_{\mathbf{u}}}
                    \mathbb{E}_{X|Y^{(i)}} X
                \right)
            -
                \mathbb{E}_{\{X^{(i)}|Y^{(i)}\}}
                f\left(
                    \frac{1}{2^{m-d}}
                    \sum_{i \in B^d_{\mathbf{u}}} X^{(i)}
                \right)
            \right.
        \right.
    \\ &\quad\quad
        \left.
        -
                f\left(
                    \frac{1}{2^{m-d}}
                    \sum_{i \in B^d_{\mathbf{u}}}
                    \mathbb{E}_{X|Y^{(\sigma(i))}} X
                \right)
            +
                \mathbb{E}_{\{X^{(i)}|Y^{(i)}\}}
                f\left(
                    \frac{1}{2^{m-d}}
                    \sum_{i \in B^d_{\mathbf{u}}} X^{(\sigma(i))}
                \right)
        \right)^2
    \\ &=
        \mathbb{E}_{\{Y^{(i)}\}}
        \left(
            \frac{1}{2^d}
            \sum_{\mathbf{u} \in \{0, 1\}^d}
            \left(
                f\left(
                    \frac{1}{2^{m-d}}
                    \sum_{i \in B^d_{\mathbf{u}}}
                    \mathbb{E}_{(X^{\prime}, X^{\prime \prime}) \sim \rho_{Y^{(i)}, Y^{(\sigma(i))}}} X^{\prime}
                \right)
            \right.
        \right.
    \\ &\quad\quad
        \left.
            \left.
            -
                \mathbb{E}_{\{(X^{\prime (i)}, X^{\prime \prime(i)}) \sim \rho_{Y^{(i)}, Y^{(\sigma(i))}}\}}
                f\left(
                    \frac{1}{2^{m-d}}
                    \sum_{i \in B^d_{\mathbf{u}}} X^{\prime (i)}
                \right)
            \right.
        \right.
    \\ &\quad\quad
        \left.
        -
            f\left(
                \frac{1}{2^{m-d}}
                \sum_{i \in B^d_{\mathbf{u}}}
                \mathbb{E}_{(X^{\prime}, X^{\prime \prime}) \sim \rho_{Y^{(i)}, Y^{(\sigma(i))}}} X^{\prime \prime}
            \right)
        \right.
    \\ &\quad\quad
        \left.
            +
            \mathbb{E}_{\{(X^{\prime (i)}, X^{\prime \prime(i)}) \sim \rho_{Y^{(i)}, Y^{(\sigma(i))}}\}}
            f\left(
                \frac{1}{2^{m-d}}
                \sum_{i \in B^d_{\mathbf{u}}} X^{\prime \prime(i)}
            \right)
        \right)^2
    \\ &\leq
        \mathbb{E}_{\{Y^{(i)}\}}
        \frac{1}{2^d}
        \sum_{\mathbf{u} \in \{0, 1\}^d}
        \mathbb{E}_{\{(X^{\prime (i)}, X^{\prime \prime(i)}) \sim \rho_{Y^{(i)}, Y^{(\sigma(i))}}\}}
    \\ &\quad\quad
        \left(
            f\left(
                \frac{1}{2^{m-d}}
                \sum_{i \in B^d_{\mathbf{u}}}
                \mathbb{E}_{(X^{\prime}, X^{\prime \prime}) \sim \rho_{Y^{(i)}, Y^{(\sigma(i))}}} X^{\prime}
            \right)
        -
            f\left(
                \frac{1}{2^{m-d}}
                \sum_{i \in B^d_{\mathbf{u}}} X^{\prime (i)}
            \right)
        \right.
    \\ &\quad\quad\quad
        \left.
        -
            f\left(
                \frac{1}{2^{m-d}}
                \sum_{i \in B^d_{\mathbf{u}}}
                \mathbb{E}_{(X^{\prime}, X^{\prime \prime}) \sim \rho_{Y^{(i)}, Y^{(\sigma(i))}}} X^{\prime \prime}
            \right)
        +
            f\left(
                \frac{1}{2^{m-d}}
                \sum_{i \in B^d_{\mathbf{u}}} X^{\prime \prime(i)}
            \right)
        \right)^2
    \\ &\leq
        \mathbb{E}_{\{Y^{(i)}\}}
        \frac{1}{2^d}
        \sum_{\mathbf{u} \in \{0, 1\}^d}
        \mathbb{E}_{\{(X^{\prime (i)}, X^{\prime \prime(i)}) \sim \rho_{Y^{(i)}, Y^{(\sigma(i))}}\}}
    \\ &\quad
        \left(
            \alpha_1
            \left\|
                \frac{1}{2^{m-d}}
                \sum_{i \in B^d_{\mathbf{u}}}
                \mathbb{E}_{(X^{\prime}, X^{\prime \prime}) \sim \rho_{Y^{(i)}, Y^{(\sigma(i))}}} X^{\prime}
            -
                \frac{1}{2^{m-d}}
                \sum_{i \in B^d_{\mathbf{u}}} X^{\prime (i)}
            \right.
        \right.
    \\ &\quad\quad
        \left.
            \left.
            -
                \frac{1}{2^{m-d}}
                \sum_{i \in B^d_{\mathbf{u}}}
                \mathbb{E}_{(X^{\prime}, X^{\prime \prime}) \sim \rho_{Y^{(i)}, Y^{(\sigma(i))}}} X^{\prime \prime}
            +
                \frac{1}{2^{m-d}}
                \sum_{i \in B^d_{\mathbf{u}}} X^{\prime \prime(i)}
            \right\|_2
        \right.
    \\ &\quad\quad
        \left.
        +
            \alpha_2
            \left\|
                \frac{1}{2^{m-d}}
                \sum_{i \in B^d_{\mathbf{u}}}
                \mathbb{E}_{(X^{\prime}, X^{\prime \prime}) \sim \rho_{Y^{(i)}, Y^{(\sigma(i))}}} X^{\prime}
            -
                \frac{1}{2^{m-d}}
                \sum_{i \in B^d_{\mathbf{u}}} X^{\prime (i)}
            \right\|_2
        \right.
    \\ &\quad\quad\quad
        \left. \times
            \left\|
                \frac{1}{2^{m-d}}
                \sum_{i \in B^d_{\mathbf{u}}}
                \mathbb{E}_{(X^{\prime}, X^{\prime \prime}) \sim \rho_{Y^{(i)}, Y^{(\sigma(i))}}} X^{\prime}
            -
                \frac{1}{2^{m-d}}
                \sum_{i \in B^d_{\mathbf{u}}}
                \mathbb{E}_{(X^{\prime}, X^{\prime \prime}) \sim \rho_{Y^{(i)}, Y^{(\sigma(i))}}} X^{\prime \prime}
            \right\|_2
        \right)^2
    \\ &\leq
        2
        \alpha_1^2
        \mathbb{E}_{\{Y^{(i)}\}}
        \frac{1}{2^d}
        \sum_{\mathbf{u} \in \{0, 1\}^d}
        \mathbb{E}_{\{(X^{\prime (i)}, X^{\prime \prime(i)}) \sim \rho_{Y^{(i)}, Y^{(\sigma(i))}}\}}
    \\ &\quad\quad
        \left\|
            \frac{1}{2^{m-d}}
            \sum_{i \in B^d_{\mathbf{u}}}
            \mathbb{E}_{(X^{\prime}, X^{\prime \prime}) \sim \rho_{Y^{(i)}, Y^{(\sigma(i))}}} X^{\prime}
        -
            \frac{1}{2^{m-d}}
            \sum_{i \in B^d_{\mathbf{u}}} X^{\prime (i)}
        \right.
    \\ &\quad\quad\quad
        \left.
        -
            \frac{1}{2^{m-d}}
            \sum_{i \in B^d_{\mathbf{u}}}
            \mathbb{E}_{(X^{\prime}, X^{\prime \prime}) \sim \rho_{Y^{(i)}, Y^{(\sigma(i))}}} X^{\prime \prime}
        +
            \frac{1}{2^{m-d}}
            \sum_{i \in B^d_{\mathbf{u}}} X^{\prime \prime (i)}
        \right\|_2^2
    \\ &\quad
    +
        2
        \alpha_2^2
        \mathbb{E}_{\{Y^{(i)}\}}
        \frac{1}{2^d}
        \sum_{\mathbf{u} \in \{0, 1\}^d}
        \mathbb{E}_{\{(X^{\prime (i)}, X^{\prime \prime(i)}) \sim \rho_{Y^{(i)}, Y^{(\sigma(i))}}\}}
    \\ &\quad\quad
        \left\|
            \frac{1}{2^{m-d}}
            \sum_{i \in B^d_{\mathbf{u}}}
            \mathbb{E}_{(X^{\prime}, X^{\prime \prime}) \sim \rho_{Y^{(i)}, Y^{(\sigma(i))}}} X^{\prime}
        -
            \frac{1}{2^{m-d}}
            \sum_{i \in B^d_{\mathbf{u}}} X^{\prime(i)}
        \right\|_2^2
    \\ &\quad\quad
        \times
        \left\|
            \frac{1}{2^{m-d}}
            \sum_{i \in B^d_{\mathbf{u}}}
            \mathbb{E}_{(X^{\prime}, X^{\prime \prime}) \sim \rho_{Y^{(i)}, Y^{(\sigma(i))}}} X^{\prime}
        -
            \frac{1}{2^{m-d}}
            \sum_{i \in B^d_{\mathbf{u}}}
            \mathbb{E}_{(X^{\prime}, X^{\prime \prime}) \sim \rho_{Y^{(i)}, Y^{(\sigma(i))}}} X^{\prime \prime}
        \right\|_2^2
    \\ &\leq
        2
        \alpha_1^2
        \mathbb{E}_{\{Y^{(i)}\}}
        \frac{1}{2^d}
        \sum_{\mathbf{u} \in \{0, 1\}^d}
        \mathbb{E}_{\{(X^{\prime (i)}, X^{\prime \prime(i)}) \sim \rho_{Y^{(i)}, Y^{(\sigma(i))}}\}}
    \\ &\quad\quad
        \left\|
            \frac{1}{2^{m-d}}
            \sum_{i \in B^d_{\mathbf{u}}}
            \mathbb{E}_{(X^{\prime}, X^{\prime \prime}) \sim \rho_{Y^{(i)}, Y^{(\sigma(i))}}} X^{\prime}
        -
            \frac{1}{2^{m-d}}
            \sum_{i \in B^d_{\mathbf{u}}} X^{\prime(i)}
        \right.
    \\ &\quad\quad
        \left.
        -
            \frac{1}{2^{m-d}}
            \sum_{i \in B^d_{\mathbf{u}}}
            \mathbb{E}_{(X^{\prime}, X^{\prime \prime}) \sim \rho_{Y^{(i)}, Y^{(\sigma(i))}}} X^{\prime \prime}
        +
            \frac{1}{2^{m-d}}
            \sum_{i \in B^d_{\mathbf{u}}} X^{\prime \prime(i)}
        \right\|_2^2
    \\ &\quad
    +
        2
        \alpha_2^2
        \mathbb{E}_{\{Y^{(i)}\}}
        \frac{1}{2^d}
        \sum_{\mathbf{u} \in \{0, 1\}^d}
        \mathbb{E}_{\{(X^{\prime (i)}, X^{\prime \prime(i)}) \sim \rho_{Y^{(i)}, Y^{(\sigma(i))}}\}}
    \\ &\quad\quad
        \left\|
            \frac{1}{2^{m-d}}
            \sum_{i \in B^d_{\mathbf{u}}}
            \mathbb{E}_{(X^{\prime}, X^{\prime \prime}) \sim \rho_{Y^{(i)}, Y^{(\sigma(i))}}} X^{\prime}
        -
            \frac{1}{2^{m-d}}
            \sum_{i \in B^d_{\mathbf{u}}} X^{\prime(i)}
        \right\|_2^2
    \\ &\quad\quad\quad
        \times
        \left(
            \frac{1}{2^{m-d}}
            \sum_{i \in B^d_{\mathbf{u}}}
            \left\|
                \mathbb{E}_{(X^{\prime}, X^{\prime \prime}) \sim \rho_{Y^{(i)}, Y^{(\sigma(i))}}} X^{\prime}
            -
                \mathbb{E}_{(X^{\prime}, X^{\prime \prime}) \sim \rho_{Y^{(i)}, Y^{(\sigma(i))}}} X^{\prime \prime}
            \right\|_2^2
        \right)
    \\ &\leq
        2
        \alpha_1^2
        \mathbb{E}_{\{Y^{(i)}\}}
        \frac{1}{2^d}
        \sum_{\mathbf{u} \in \{0, 1\}^d}
        \frac{2}{2^{2(m-d)}} \sum_{i \in B^d_{\mathbf{u}}}
        \mathbb{E}_{\{(X^{\prime (i)}, X^{\prime \prime(i)}) \sim \rho_{Y^{(i)}, Y^{(\sigma(i))}}\}}
        \left\|
            X^{\prime (i)} - X^{\prime \prime(i)}
        \right\|_2^2
    \\ &\quad\quad
    +
        2
        \alpha_2^2
        \mathbb{E}_{\{Y^{(i)}\}}
        \frac{1}{2^d}
        \sum_{\mathbf{u} \in \{0, 1\}^d}
        \frac{1}{2^{2(m-d)}}
        \sum_{i \in B^d_{\mathbf{u}}}
        \mathbb{V}_{X|Y^{(i)}} X
    \\ &\quad\quad\quad \times
        \left(
            \frac{1}{2^{m-d}}
            \sum_{i \in B^d_{\mathbf{u}}}
            \mathbb{E}_{(X^{\prime}, X^{\prime \prime}) \sim \rho_{Y^{(i)}, Y^{(\sigma(i))}}}
            \left\|
                X^{\prime} - X^{\prime \prime}
            \right\|_2^2
        \right)
    \\ &\leq
        2
        \alpha_1^2
        \beta^2
        \mathbb{E}_{\{Y^{(i)}\}}
        \frac{1}{2^d}
        \sum_{\mathbf{u} \in \{0, 1\}^d}
        \frac{1}{2^{2(m-d)}}
        \sum_{i \in B^d_{\mathbf{u}}}
        \left\|
            \mathcal{T}(Y^{(i)})
        -
            \mathcal{T}(Y^{(\sigma(i))})
        \right\|^2_2
    \\ &\quad\quad
    +
        2
        \alpha_2^2
        \beta^2
        \mathbb{E}_{\{Y^{(i)}\}}
        \frac{1}{2^d}
        \sum_{\mathbf{u} \in \{0, 1\}^d}
        \frac{1}{2^{2(m-d)}}
        \sum_{i \in B^d_{\mathbf{u}}}
        \mathbb{V}_{X|Y^{(i)}} X
        \times
        \left(
        \frac{1}{2^{m-d}}
        \sum_{i \in B^d_{\mathbf{u}}}
        \left\|
            \mathcal{T}(Y^{(i)})
        -
            \mathcal{T}(Y^{(\sigma(i))})
        \right\|_2^2
        \right)
    \\ &\leq
        2\beta^2
        \left(
            \alpha_1^2
        +
            \alpha_2^2
            \sup_{Y}\mathbb{V}_{X|Y} X
        \right)
        \frac{1}{2^{m-d}}
        \mathbb{E}_{\{Y^{(i)}\}}
        \frac{1}{2^d}
        \sum_{\mathbf{u} \in \{0, 1\}^d}
        \frac{1}{2^{m-d}}
        \sum_{i \in B^d_{\mathbf{u}}}
        \left\|
            \mathcal{T}(Y^{(i)})
        -
            \mathcal{T}(Y^{(\sigma(i))})
        \right\|_2^2
    \\ &\leq
        2\beta^2
        \left(
            \alpha_1^2
        +
            \alpha_2^2
            \sup_{Y}\mathbb{V}_{X|Y} X
        \right)
        \frac{1}{2^{m-d}}
        \mathbb{E}_{\{Y^{(i)}\}}
        \frac{1}{2^{d-1}}
        \sum_{\mathbf{u} \in \{0, 1\}^{d-1}}
        \max_{i, j \in B^{d-1}_{\mathbf{u}}}
        \left\|
            \mathcal{T}(Y^{(i)})
        -
            \mathcal{T}(Y^{(j)})
        \right\|_2^2
    \\ &\leq
        2\beta^2
        \left(
            \alpha_1^2
        +
            \alpha_2^2
            \sup_{Y}\mathbb{V}_{X|Y} X
        \right)
        \frac{1}{2^{m-d}}\cdot
        \frac{2K}{2^{(d-1)/K}}
        \leq
        2\beta^2
        \left(
            \alpha_1^2
        +
            \alpha_2^2
            \sup_{Y}\mathbb{V}_{X|Y} X
        \right)
        \frac{2K}{2^{(m-1)/K}}.
    \end{align}

    Therefore, the second term of \eqref{eq:main_theorem_three_terms} is bounded above by
    \begin{align}
    &\quad
        \mathbb{E}_{\{Y^{(i)}\}}
        \left(
            \left(
                \sum_{d=0}^{m}
                \hat{P}_{\mathrm{in}}^{(d,m-d)}
            -
                \sum_{d=1}^{m}
                \hat{Q}_{\mathrm{in}}^{(d,m-d)}
            \right)
        -
            \left(
                \sum_{d=0}^{m}
                \hat{P}_{\mathrm{out}}^{(d,m-d)}
            -
                \sum_{d=1}^{m}
                \hat{Q}_{\mathrm{out}}^{(d,m-d)}
            \right)
        \right)^2
    \\ &\leq
       (m+1) \mathbb{E}_{\{Y^{(i)}\}}
        \left(
            \hat{P}_{\mathrm{in}}^{(0,m)}-\hat{P}_{\mathrm{out}}^{(0,m)}
        \right)^2
    \\ &\quad\quad
        +
       (m+1) \sum_{d=1}^{m}\mathbb{E}_{\{Y^{(i)}\}}
        \left(
                \hat{P}_{\mathrm{in}}^{(d,m-d)}
            -
                \hat{P}_{\mathrm{out}}^{(d,m-d)}
        -
                \hat{Q}_{\mathrm{in}}^{(d,m-d)}
            +
                \hat{Q}_{\mathrm{out}}^{(d,m-d)}
        \right)^2
    \\ & \leq
        (m+1)\frac{\alpha_1^2}{2^m}
        \mathbb{E}_{Y}\mathbb{V}_{X|Y} X
        +(m+1)\sum_{d=0}^{m-1}2\beta^2
        \left(
            \alpha_1^2
        +
            \alpha_2^2
            \sup_{Y}\mathbb{V}_{X|Y} X
        \right)
        \frac{2K}{2^{(m-1)/K}}
        =
        O\left(
            \frac{m^2}{2^{m/K}}
        \right).
    \end{align}

\underline{A bound on the third term of \eqref{eq:main_theorem_three_terms}:}
    Finally a bound on the third term of \eqref{eq:main_theorem_three_terms} is bounded from above as follows.
    \begin{align}
    &\quad
        \mathbb{E}_{\{(X^{(i)},Y^{(i)})\}}
        \left(
            \left(
                \sum_{d=0}^{m}
                \hat{P}_{\mathrm{out}}^{(d,m-d)}
            -
                \sum_{d=1}^{m}
                \hat{Q}_{\mathrm{out}}^{(d,m-d)}
            \right)
        -
            \left(
                \sum_{d=0}^{m}
                \hat{P}^{(d,m-d)}
            -
                \sum_{d=1}^{m}
                \hat{Q}^{(d,m-d)}
            \right)
        \right)^2
    \\ &\leq
        (2m+1)\sum_{d=0}^{m}\mathbb{E}_{\{(X^{(i)},Y^{(i)})\}}\left( \hat{P}_{\mathrm{out}}^{(d,m-d)}-\hat{P}^{(d,m-d)}\right)^2
    \\ &\quad\quad
        + (2m+1)\sum_{d=1}^{m}\mathbb{E}_{\{(X^{(i)},Y^{(i)})\}}\left( \hat{Q}_{\mathrm{out}}^{(d,m-d)}-\hat{Q}^{(d,m-d)}\right)^2
    \\ &=
        (2m+1)
        \sum_{d=0}^{m}
        \mathbb{E}_{\{(X^{(i)},Y^{(i)})\}}
    \\ &\quad\quad
        \left(
            \frac{1}{2^d}
            \sum_{\mathbf{u} \in \{0, 1\}^d}
            \left(
                \mathbb{E}_{\{X^{(i)}|Y^{(i)}\}}
                f\left(
                    \frac{1}{2^{m-d}}
                    \sum_{i \in B^d_{\mathbf{u}}} X^{(i)}
                \right)
            -
                f\left(
                    \frac{1}{2^{m-d}}
                    \sum_{i \in B^d_{\mathbf{u}}} X^{(i)}
                \right)
            \right)
        \right)^2
    \\ &\quad\quad
    +
        (2m+1)
        \sum_{d=1}^{m}
        \mathbb{E}_{\{(X^{(i)},Y^{(i)})\}}
    \\ &\quad\quad\quad
        \left(
            \frac{1}{2^d}
            \sum_{\mathbf{u} \in \{0, 1\}^d}
            \left(
                \mathbb{E}_{\{X^{(i)}|Y^{(i)}\}}
                f\left(
                    \frac{1}{2^{m-d}}
                    \sum_{i \in B'^d_{\mathbf{u}}} X^{(i)}
                \right)
            -
                f\left(
                    \frac{1}{2^{m-d}}
                    \sum_{i \in B'^d_{\mathbf{u}}} X^{(i)}
                \right)
            \right)
        \right)^2
    \\ &=
        (2m+1)
        \sum_{d=0}^{m}
        \mathbb{E}_{\{Y^{(i)}\}}
        \frac{1}{2^{2d}}
        \sum_{\mathbf{u} \in \{0, 1\}^d}
        \mathbb{V}_{\{X^{(i)}|Y^{(i)}\}}
        f\left(
            \frac{1}{2^{m-d}}
            \sum_{i \in B^d_{\mathbf{u}}} X^{(i)}
        \right)
    \\ &\quad
    +
        (2m+1)
        \sum_{d=1}^{m}
        \mathbb{E}_{\{Y^{(i)}\}}
        \frac{1}{2^{2d}}
        \sum_{\mathbf{u} \in \{0, 1\}^d}
        \mathbb{V}_{\{X^{(i)}|Y^{(i)}\}}
        f\left(
            \frac{1}{2^{m-d}}
            \sum_{i \in B'^d_{\mathbf{u}}} X^{(i)}
        \right)
    \\ &\leq
        \alpha_1^2
        (2m+1)
        \sum_{d=0}^{m}
        \mathbb{E}_{\{Y^{(i)}\}}
        \frac{1}{2^{2d}}
        \sum_{\mathbf{u} \in \{0, 1\}^d}
        \mathbb{V}_{\{X^{(i)}|Y^{(i)}\}}
        \frac{1}{2^{m-d}}
        \sum_{i \in B^d_{\mathbf{u}}} X^{(i)}
    \\ &\quad
    +
        \alpha_1^2
        (2m+1)
        \sum_{d=1}^{m}
        \mathbb{E}_{\{Y^{(i)}\}}
        \frac{1}{2^{2d}}
        \sum_{\mathbf{u} \in \{0, 1\}^d}
        \mathbb{V}_{\{X^{(i)}|Y^{(i)}\}}
        \frac{1}{2^{m-d}}
        \sum_{i \in B'^d_{\mathbf{u}}} X^{(i)}
    \\ &=
        \alpha_1^2
        (2m+1)
        \sum_{d=0}^{m}
        \mathbb{E}_{\{Y^{(i)}\}}
        \frac{1}{2^{2d}}
        \sum_{\mathbf{u} \in \{0, 1\}^d}
        \frac{1}{2^{2(m-d)}}
        \sum_{i \in B^d_{\mathbf{u}}}
        \mathbb{V}_{\{X^{(i)}|Y^{(i)}\}} X^{(i)}
    \\ &\quad
    +
        \alpha_1^2
        (2m+1)
        \sum_{d=1}^{m}
        \mathbb{E}_{\{Y^{(i)}\}}
        \frac{1}{2^{2d}}
        \sum_{\mathbf{u} \in \{0, 1\}^d}
        \frac{1}{2^{2(m-d)}}
        \sum_{i \in B'^d_{\mathbf{u}}}
        \mathbb{V}_{\{X^{(i)}|Y^{(i)}\}} X^{(i)}
    \\ &=
        \alpha_1^2
        (2m+1)
        \sum_{d=0}^{m}
        \mathbb{E}_{\{Y^{(i)}\}}
        \frac{1}{2^{2m}}
        \sum_{i=1}^{2^m}
        \mathbb{V}_{\{X^{(i)}|Y^{(i)}\}} X^{(i)}
    \\ &\quad
    +
        \alpha_1^2
        (2m+1)
        \sum_{d=1}^{m}
        \mathbb{E}_{\{Y^{(i)}\}}
        \frac{1}{2^{2m}}
        \sum_{i=1}^{2^m}
        \mathbb{V}_{\{X^{(i)}|Y^{(i)}\}} X^{(i)}
    \\ &=
        \alpha_1^2
        (2m+1)^2
        \mathbb{E}_{\{Y^{(i)}\}}
        \frac{1}{2^{2m}}
        \sum_{i=1}^{2^m}
        \mathbb{V}_{\{X^{(i)}|Y^{(i)}\}} X^{(i)}
        =
        \frac{1}{2^m}
        \alpha_1^2
        (2m+1)^2\,
        \mathbb{E}_{Y}
        \mathbb{V}_{X|Y} X
        = O\left( \frac{m^2}{2^m}\right).
    \end{align}
    Combining the respective upper bound on each of the three terms of \eqref{eq:main_theorem_three_terms}, we complete the proof of the theorem.
\end{proof}

Our theoretical result shows that the MSE of our proposed method decays at least as fast as $N^{-1/K}(\log N)^2$. Although we omit the details, in comparison, the authors' previous method in \cite{hironaka2023efficient} has an upper bound of $O(N^{-1/2K})$ under the same conditions as in Theorem~\ref{thm:main}. Thus, our proposed method has a faster rate of convergence in terms of upper bounds. The rate $N^{-1/K}(\log N)^2$ indicates that our method is efficient when the dimensionality of the outer random variable $K$ is small. However, due to the exponential dependence on $K$, our method may not be efficiently applicable to high-dimensional problems, and estimating nested expectations may suffer from the curse of dimensionality.

In fact, a similar exponential dependence of the approximation error on $K$ also appears in other studies, such as \cite{hong2017kernel} and \cite{broto2020variance}, in which nested expectations are estimated without sampling from the inner conditional distribution. This implies that the curse of dimensionality is a common issue in estimating nested expectations without inner conditional sampling. Further studies could investigate whether it is possible to overcome the curse of dimensionality or not.

%% file: 40-numerical.tex
\section{Numerical Experiments}\label{sec:numerical}

Here we present numerical experiments to evaluate the performance of our proposed method and compare it with some existing methods. As mentioned in Section~\ref{sec:intro} as a motivating example for nested expectations, we focus on the problem of estimating EVSI. In what follows, we first introduce two test problems in Section~\ref{ssec:seting}, and then report the numerical results and make some discussion in Section~\ref{ssec:result}.

\input{42-settings}

\input{43-results}

\input{50-discussion}

%% file: 42-settings.tex
\subsection{Test Problem}\label{ssec:seting}

The first example is somewhat artificial, but intentionally designed to be simple enough to test whether the estimates from our proposed method converge to the correct value. In what follows, we denote the Bernoulli distribution with success probability $p$ by $\mathrm{B}(1, p)$.

\begin{problem}[simple test case]\label{prm:evsi0}
Consider model parameters $p\in (0,1)$ and $M\in \mathbb{N}$, and let $\theta \sim \mathrm{B}\left(1, \frac{1}{2}\right)$ be an inner random variable. Define the set $D=\{0,1\}$ and, for each option $d\in D$, let $\mathrm{NB}_d$ be a net benefit function given by
\begin{align}
\mathrm{NB}_d = \begin{cases}
\theta, & \text{if $d = 0$,} \\
1 - \theta, & \text{if $d = 1$.}
\end{cases}
\end{align}
We define the sample information $Y$ as an $M$-dimensional vector $(Y_1, \dots, Y_M)$, where
\begin{align}
Y_m = (2b_m-1) U_m (2\theta-1),
\end{align}
for $m \in \{1, \dots, M\}$. Here, $\{b_m\}_m$ and $\{U_m\}_m$ are sets of independent random variables that follow the Bernoulli distribution $\mathrm{B}\left(1, p\right)$ and the uniform distribution $\mathrm{U}\left(0, 1\right)$, respectively.
\end{problem}

For this problem, the exact value of EVSI is available, which is given by
\[ I=\frac{1}{2}\sum_{m=0}^M\binom{M}{m}\max\left(p^{M-m}(1-p)^M, p^M(1-p)^{M-m}\right).\]
Using this, we can estimate the MSE of the estimator $\hat{I}$ without bias as
\begin{align}\label{eq:estimate_MSE}
    \frac{1}{r}\sum_{i=1}^{r}\left( I-\hat{I}^{(i)}\right)^2,
\end{align}
where $r\in \mathbb{N}$ and $\hat{I}^{(1)},\dots,\hat{I}^{(r)}$ denote $r$ independent estimates obtained by the estimator $\hat{I}$. In what follows, we always choose $r=100$.

For the second example, let us consider a medical decision model presented in \cite{reeves2019three}. This model concerns the determination of the optimal dressing types (or no dressing) to reduce surgical site infections (SSI) in primary surgical wounds and aims to assess the feasibility of a multicenter randomized controlled trial (RCT). We aim to demonstrate the effectiveness of our proposed method in practical settings by applying it to this realistic decision-making problem.

\begin{problem}[medical decision problem from \cite{reeves2019three}]\label{prm:evsi2}
    Let $\theta$ be a vector of inner random variables, the description and probability distribution of each of which are described in Table~\ref{tbl:bluebelle_theta}.
    \begin{table}[t]
    \caption{Model inputs for the second problem presented in \cite{reeves2019three}. Note that $\text{lognormal}(\mu,\Sigma)$ represents the log-normal distribution, where $\mu$ and $\Sigma$ denote the mean vector and covariance matrix of the corresponding normal distribution, respectively.}
    \label{tbl:bluebelle_theta}
    \centering
    \small
    \begin{tabular}{|l|l|l|}
            \hline
            Description & Parameter & Distribution \\
            \hline \hline
            Willingness to pay & $\mathit{WTP}$ & $20000$ (constant) \\
            per QALY threshold & & \\ \hline
            Quality adjusted life-years & $\mathit{SSIQALYloss}$ & $0.12$ (constant) \\
            decrement resulting from an SSI & & \\ \hline
            Cost attributable to SSI & $\mathit{SSIcost}$ & $\mathrm{lognormal}\left(8.972, 0.1631^2 \right)$ \\ \hline
            Cost for no dressing & $\mathit{dressingcosts}_E$ & $0$ (constant) \\ \hline
            Cost for simple dressings & $\mathit{dressingcosts}_S$ & $5.25$ (constant) \\ \hline
            Cost for glue dressing & $\mathit{dressingcosts}_G$ & $13.86$ (constant) \\ \hline
            Cost for advanced dressings & $\mathit{dressingcosts}_A$ & $21.39$ (constant) \\ \hline
            SSI risk with simple dressings & $\mathit{pSSI}_S$ & $\mathrm{N}\left(0.1380, 0.0018^2 \right)$ \\ \hline
            SSI risk with other dressing & $\mathit{pSSI}_d$ & Derived from $\mathit{pSSI}_S$ and $\mathit{OR}_d$ \\
            types $(d\in \{E,G,A\})$ & & \\ \hline
            Odds ratios of SSI risk & \multirow{3}{*}{$\begin{pmatrix} \mathit{OR}_E \\ \mathit{OR}_G \\ \mathit{OR}_A \end{pmatrix}$} & \multirow{3}{*}{$\mathrm{lognormal}\left(
                \begin{pmatrix} -0.05 \\ -0.07 \\ -0.18 \end{pmatrix},
                \begin{pmatrix}
                    0.07 & 0.06 & 0.02 \\
                    0.06 & 0.22 & 0.02 \\
                    0.02 & 0.02 & 0.05
                \end{pmatrix}
            \right) $} \\
            relative to simple dressings & & \\
            & & \\ \hline
        \end{tabular}
    \end{table}
    Consider a set of four different dressing types, denoted by $D = \{E, S, G, A\}$, and for each type $d \in D$, we define the net benefit function as follows:
    \begin{align}
        \mathrm{NB}_d = -\left(
            \mathit{dressingcost}_d
            + \mathit{pSSI}_d
            * \left(
                \mathit{SSIcost}
                + \mathit{SSIQALYloss} \times \mathit{WTP}
            \right)
        \right).
    \end{align}
    Consider an RCT with $n$ participants, in which each treatment $d\in D$ is applied to $n_d$ participants, where $n=n_E+n_S+n_G+n_A$. We model the sample information $Y$ obtained by the RCT as a $3$-dimensional vector where
    \begin{align}
        \begin{pmatrix} Y_0 \\ Y_1 \\ Y_2 \end{pmatrix}
        \sim
        \mathrm{Normal}\left(
            \begin{pmatrix} \log(\mathit{OR}_E) \\ \log(\mathit{OR}_G) \\ \log(\mathit{OR}_A) \end{pmatrix},
            \begin{pmatrix}
                \frac{s^2 (n_S + n_E)}{n_S n_E} & \frac{s^2}{n_S} & \frac{s^2}{n_S} \\
                \frac{s^2}{n_S} & \frac{s^2 (n_S + n_G)}{n_S n_G} & \frac{s^2}{n_S} \\
                \frac{s^2}{n_S} & \frac{s^2}{n_S} & \frac{s^2 (n_S + n_A)}{n_S n_A}
            \end{pmatrix}
        \right).
    \end{align}
    where $s = 3.7$ is the standard deviation for the log-odds ratios on a given arm.
\end{problem}

In our experiments, we consider the following two scenarios
\begin{itemize}
    \item \textbf{EvSvG:} $n_E = n_S = 2n/5$, $n_G = n/5$ and $n_A=0$    
    \item \textbf{EvSvGvA:} $n_E = n_S = n_G = n_A = n/4$
\end{itemize}
with various values of $n$. In cases where $n_A=0$, we substitute a sufficiently small value of $10^{-3}$ for $n_A$ to avoid division by zero errors during computation. Although the exact value of EVSI cannot be calculated analytically for this model, we estimate it with the standard error around $0.05$ using multilevel Monte Carlo methods \cite{giles2019decision,hironaka2020multilevel} as a reference. Denoting this estimated value by $I$, we estimate the MSE of each estimator $\hat{I}$ in the same way as shown in \eqref{eq:estimate_MSE} with $r=100$.

It should be noted that the components of the outer random variables are not mutually independent. While the previous method proposed by the authors \cite{hironaka2023efficient} did not account for such situations in its theoretical analysis, the analysis presented in Section~\ref{sec:theoretical} provides a bound on the MSE even when the outer variables are dependent. Therefore, we expect that the method proposed in this paper will perform better than the previous one. Additionally, although i.i.d.\ samples from the inner conditional distribution $\rho(\theta|Y)$ can be easily generated in these problems, this assumption does not hold in general problem settings. 

In this paper, we compare our proposed method with two existing methods: the regression-based method proposed in \cite{strong2014estimating,strong2015estimating}, and the previous method proposed by the authors \cite{hironaka2023efficient}. Specifically, we use the generalized additive model (GAM) for the regression-based method and refer to it as the GAM-based method. We confirmed beforehand that using Gaussian process regression for the regression-based method would be too computationally expensive to complete the experiments within a reasonable time frame. Regarding the previous method by the authors, we consider its slightly modified version, i.e., $\hat{P}^{(m/2,m/2)}$, as we stated in Section~2, which enables us a fair comparison in terms of the total number of samples. The source codes used in our experiments are available at \url{https://github.com/Goda-Research-Group/experiment-2023-06}.

%% file: 43-results.tex
\subsection{Results and discussion}\label{ssec:result}

\begin{figure}
     \centering
     \begin{subfigure}[b]{0.49\textwidth}
         \centering
         \includegraphics[width=\textwidth]{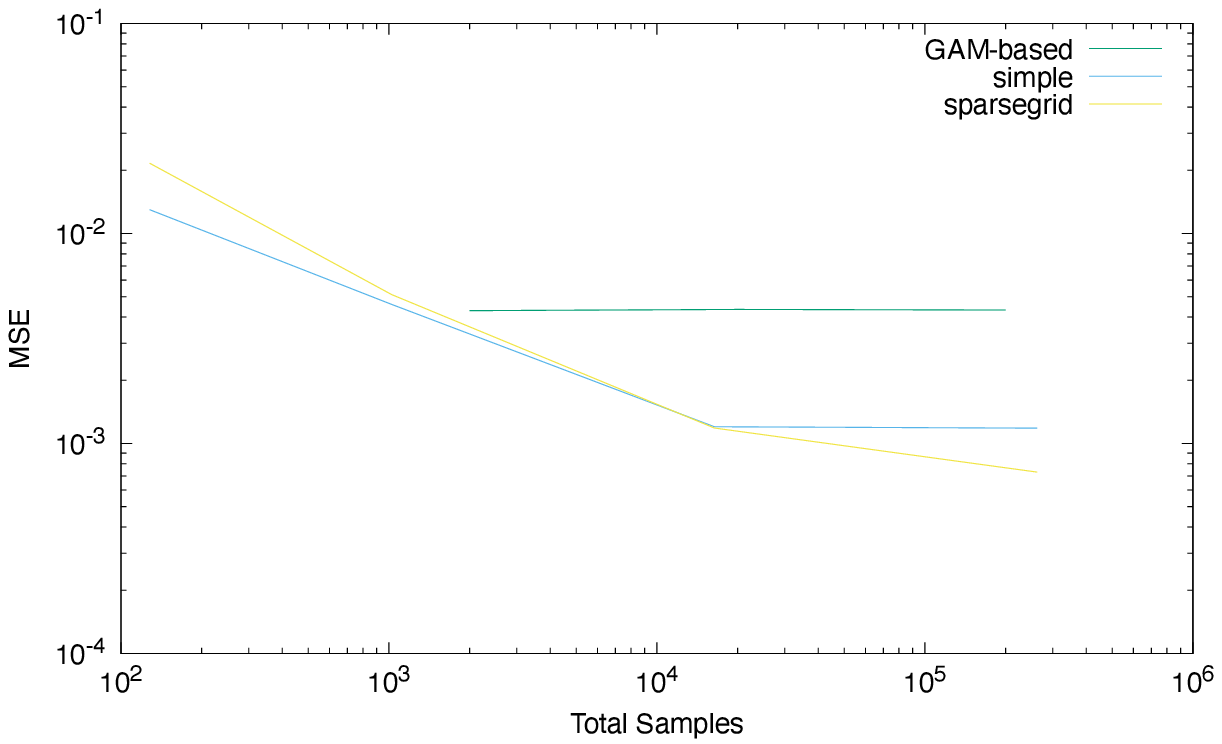}
         \caption{$(M,p)=(7,0.7)$}
     \end{subfigure}
     \hfill
     \begin{subfigure}[b]{0.49\textwidth}
         \centering
         \includegraphics[width=\textwidth]{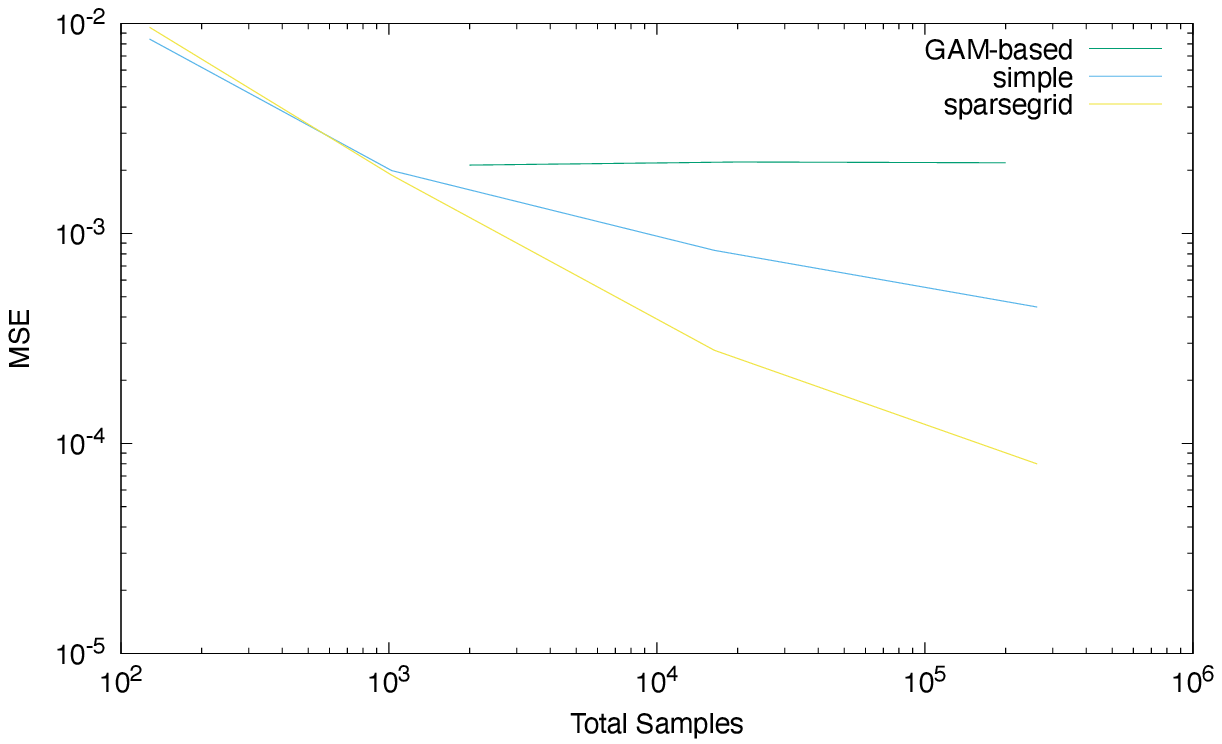}
         \caption{$(M,p)=(7,0.9)$}
     \end{subfigure}
     \\
     \begin{subfigure}[b]{0.49\textwidth}
         \centering
         \includegraphics[width=\textwidth]{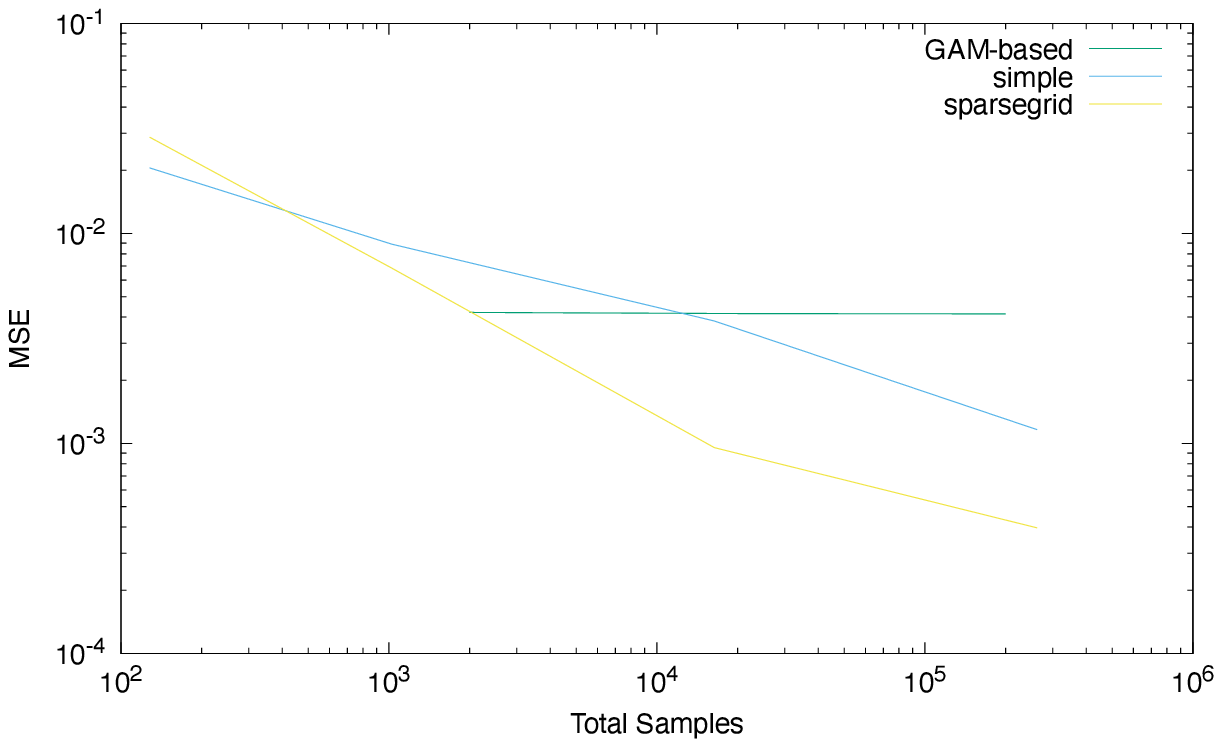}
         \caption{$(M,p)=(10,0.7)$}
     \end{subfigure}
     \hfill
     \begin{subfigure}[b]{0.49\textwidth}
         \centering
         \includegraphics[width=\textwidth]{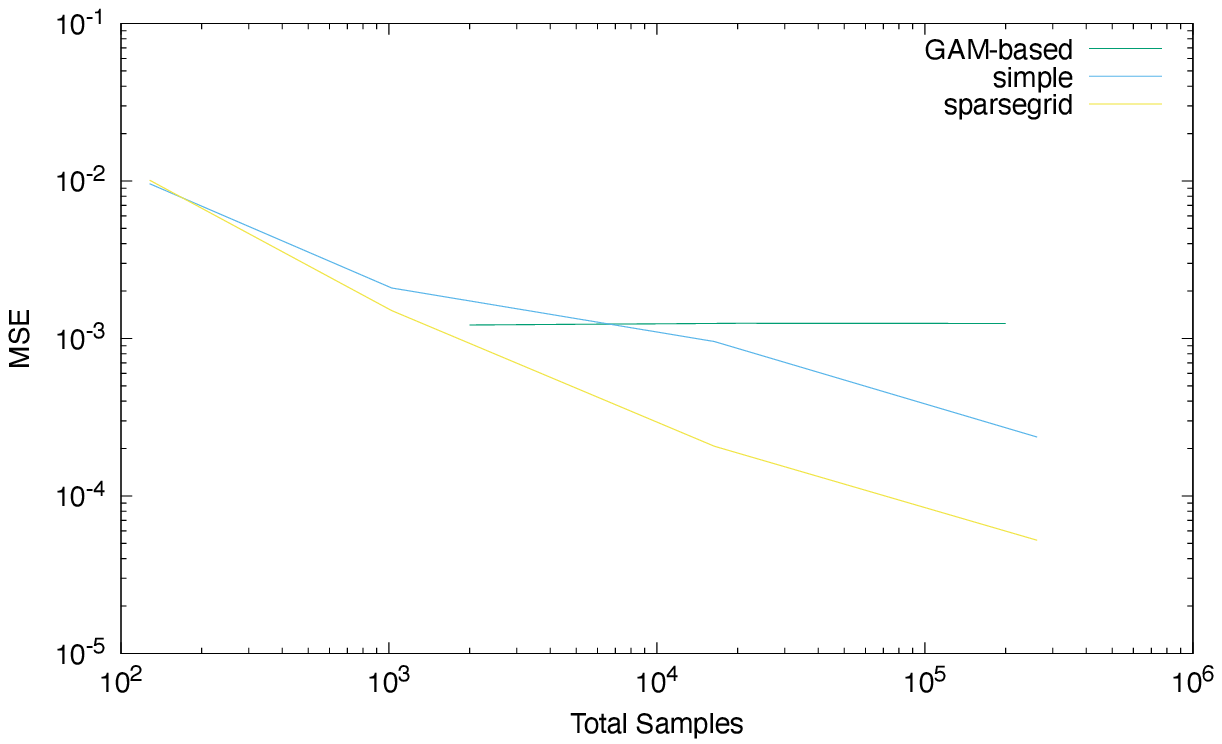}
         \caption{$(M,p)=(10,0.9)$}
     \end{subfigure}
    \caption{Comparison of MSE between the three methods (the GAM-based method, our previous and proposed methods) for the first test problem with various $(M,p)$'s. Our previous and proposed methods are denoted by {\tt simple} and {\tt sparsegrid}, respectively.}
    \label{fig:result_first}
\end{figure}

In the first experiment, we estimate the EVSI of the first model (Problem~\ref{prm:evsi0}) with the model parameters $(M,p) = (7, 0.7), (7, 0.9), (10, 0.7)$ and $(10, 0.9)$. The results are shown in the subplots of Figure~\ref{fig:result_first}, respectively. In each subplot, the vertical axis represents the MSE, and the horizontal axis represents the total number of samples used for the pair of inner and outer random variables $(\theta,Y)$.  These results demonstrate that our two methods converge to the correct EVSI value, whereas the MSE for the GAM-based method does not decay towards 0 because the estimate converges to a wrong value. Comparing our two methods, we can see that our proposed method achieves a lower error than the previous method, which is consistent with what our theoretical result predicts. Although $M$ denotes the dimensionality of the outer random variables, the convergence behavior of our proposed method is better than what is expected from the theoretical result.

\begin{figure}
     \centering
     \begin{subfigure}[b]{0.49\textwidth}
         \centering
         \includegraphics[width=\textwidth]{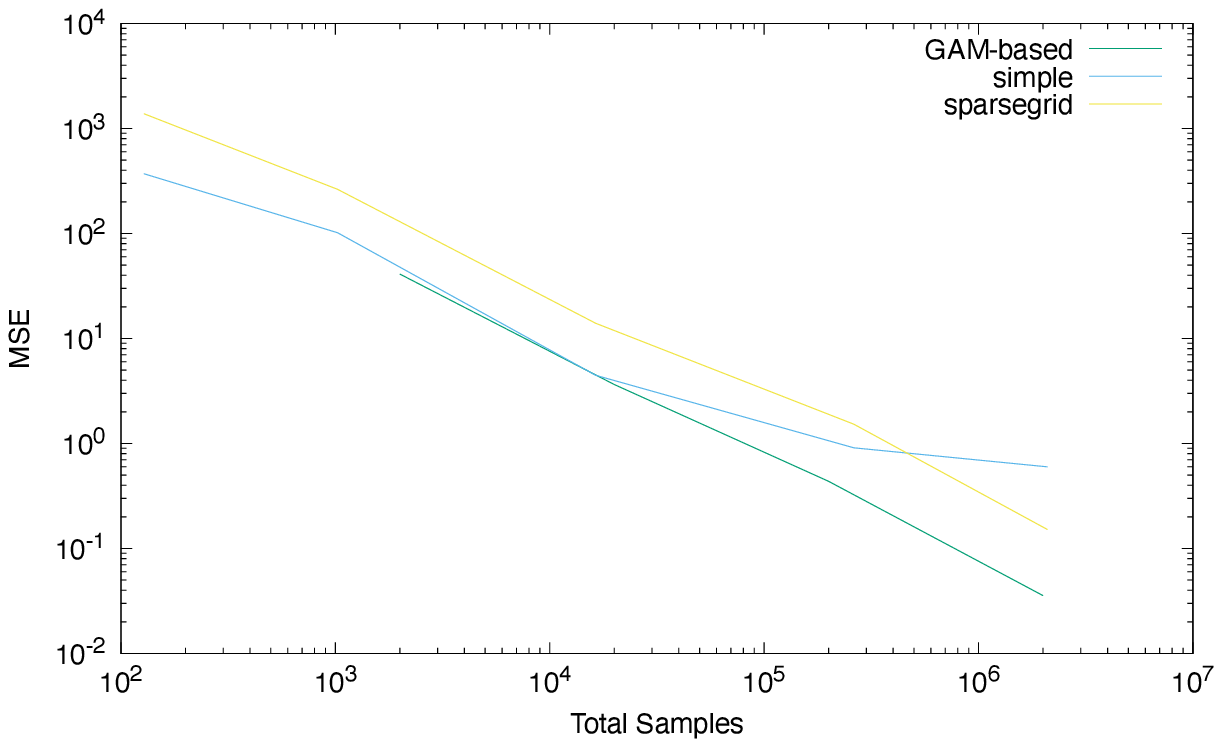}
         \caption{$n=250$}
     \end{subfigure}
     \hfill
     \begin{subfigure}[b]{0.49\textwidth}
         \centering
         \includegraphics[width=\textwidth]{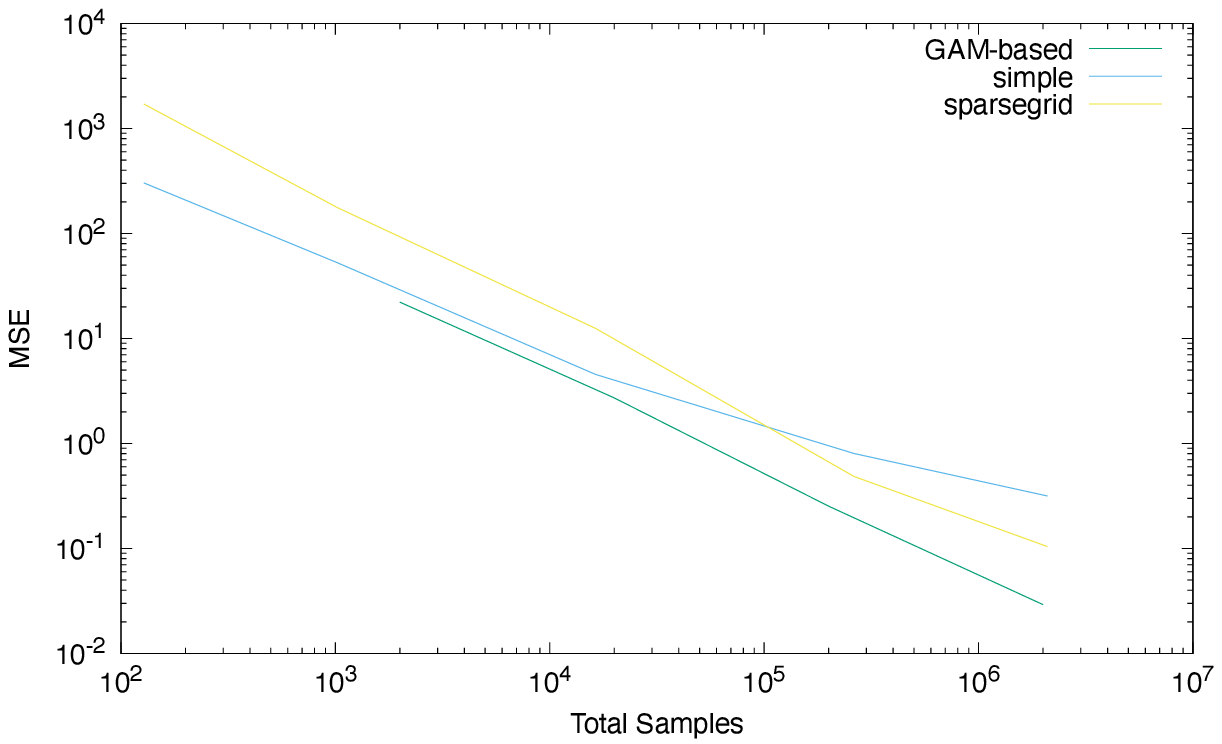}
         \caption{$n=500$}
     \end{subfigure}
     \\
     \begin{subfigure}[b]{0.49\textwidth}
         \centering
         \includegraphics[width=\textwidth]{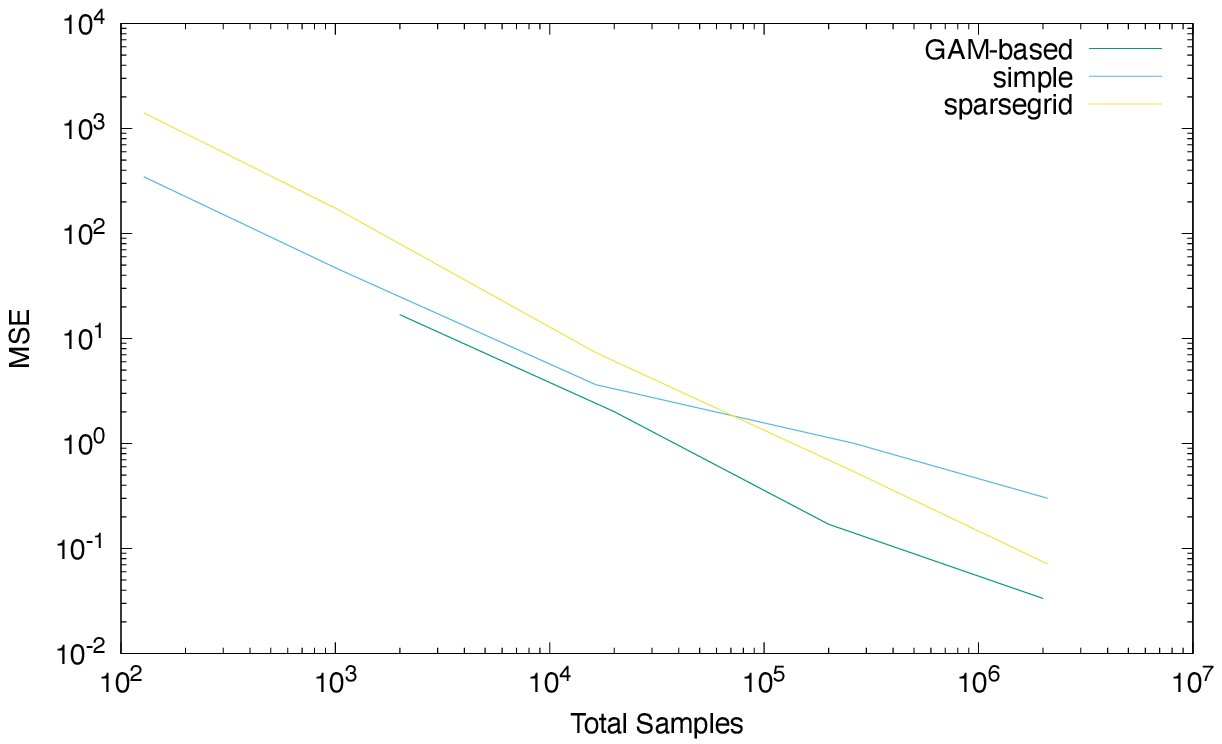}
         \caption{$n=1000$}
     \end{subfigure}
     \hfill
     \begin{subfigure}[b]{0.49\textwidth}
         \centering
         \includegraphics[width=\textwidth]{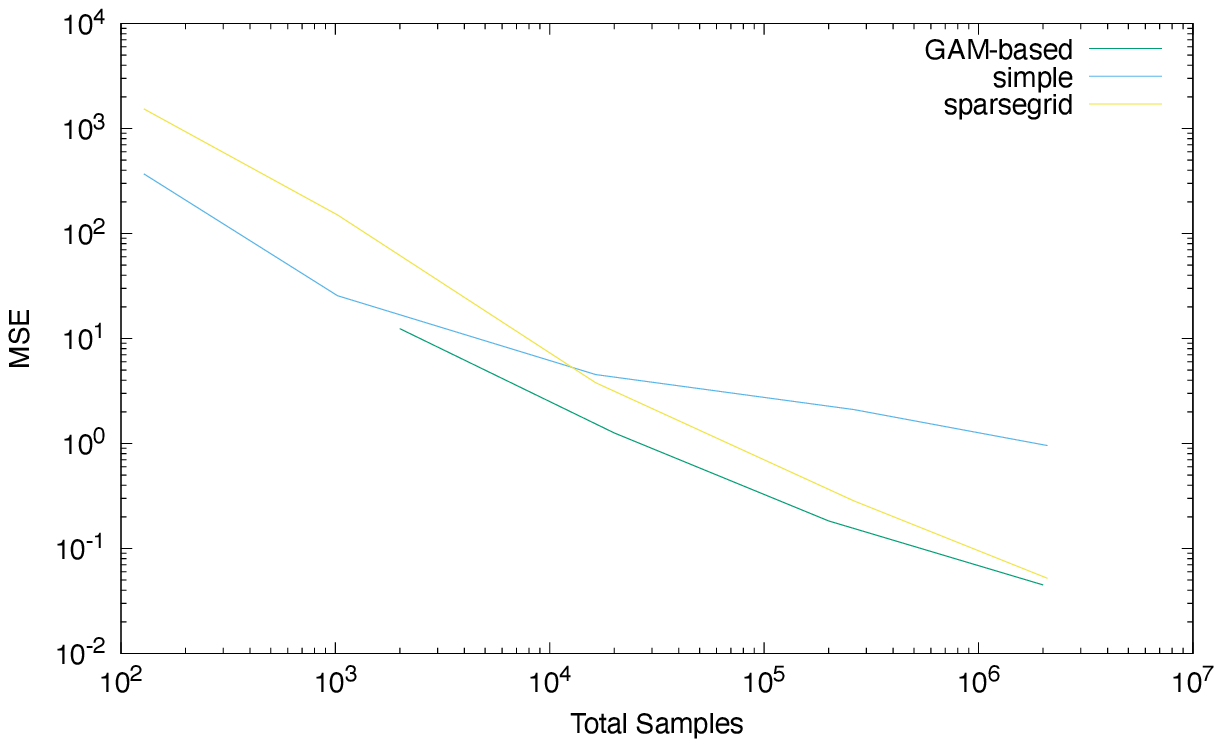}
         \caption{$n=2000$}
     \end{subfigure}
    \caption{Comparison of MSE between the three methods (the GAM-based method, our previous and proposed methods) for the scenario EvSvG with various values of $n$. Our previous and proposed methods are denoted by {\tt simple} and {\tt sparsegrid}, respectively.}
    \label{fig:result_esg221}
\end{figure}

We now move on to the second example. The results of estimating EVSI for the scenario EvSvG with the total number of participants $n=200, 500, 1000, 2000$ are presented in the subplots of Figure~\ref{fig:result_esg221}, respectively. In this scenario, it can be confirmed that the GAM-based method exhibits the best performance in terms of MSE with respect to the number of total samples for any value of $n$. However, the MSE of our proposed method also converges to 0 at a rate similar to that of the GAM-based method, making it competitive. In contrast, our previous method exhibits a lower convergence rate and clear underperformance, especially in the regime of large number of total samples.

\begin{figure}
     \centering
     \begin{subfigure}[b]{0.49\textwidth}
         \centering
         \includegraphics[width=\textwidth]{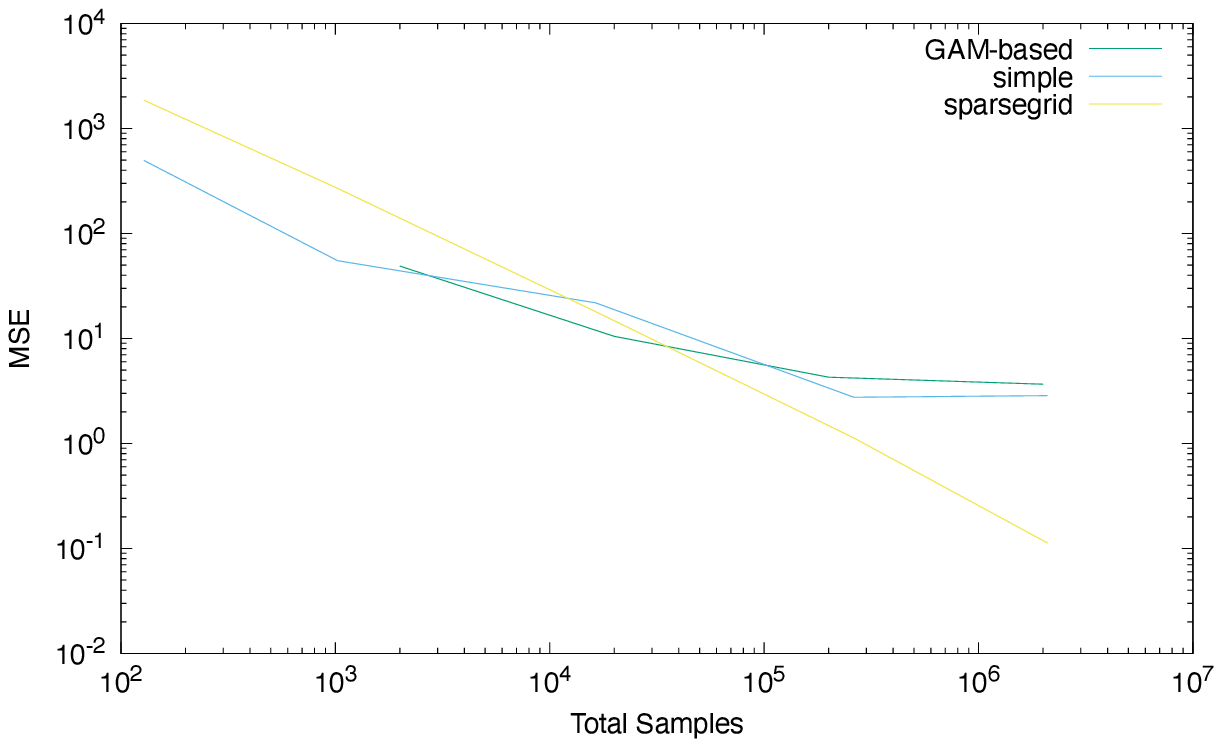}
         \caption{$n=200$}
     \end{subfigure}
     \hfill
     \begin{subfigure}[b]{0.49\textwidth}
         \centering
         \includegraphics[width=\textwidth]{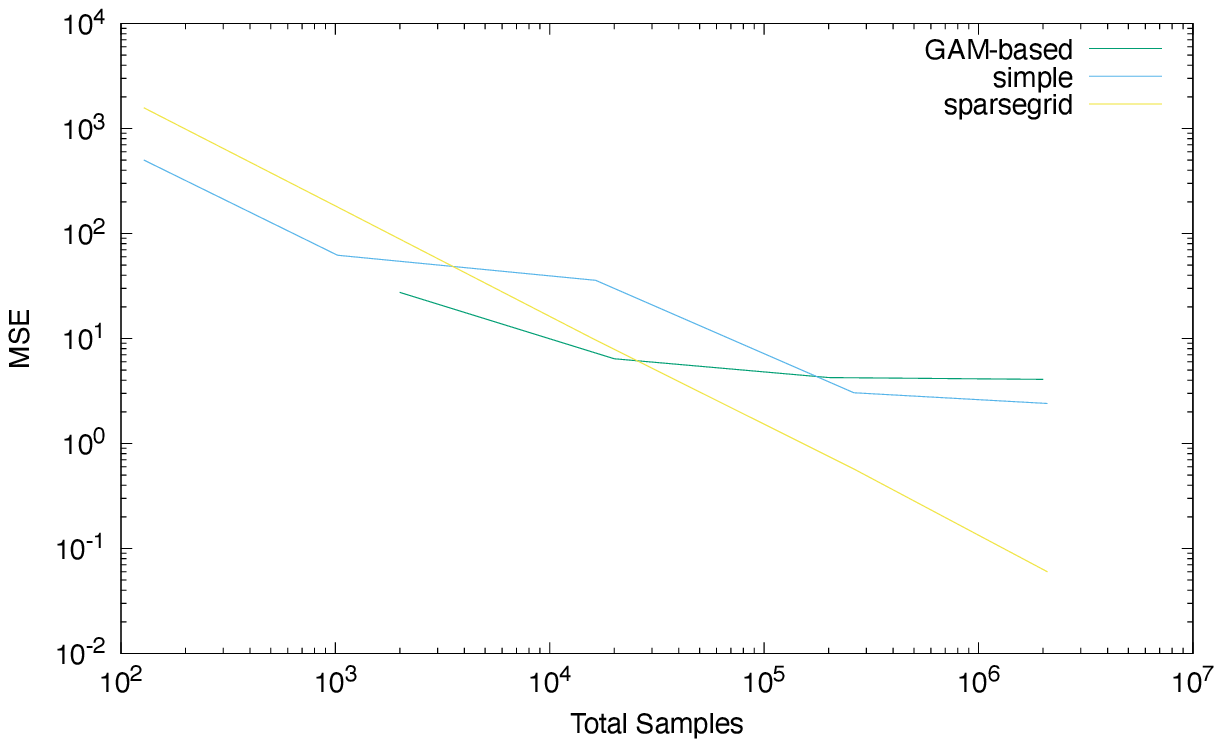}
         \caption{$n=500$}
     \end{subfigure}
     \\
     \begin{subfigure}[b]{0.49\textwidth}
         \centering
         \includegraphics[width=\textwidth]{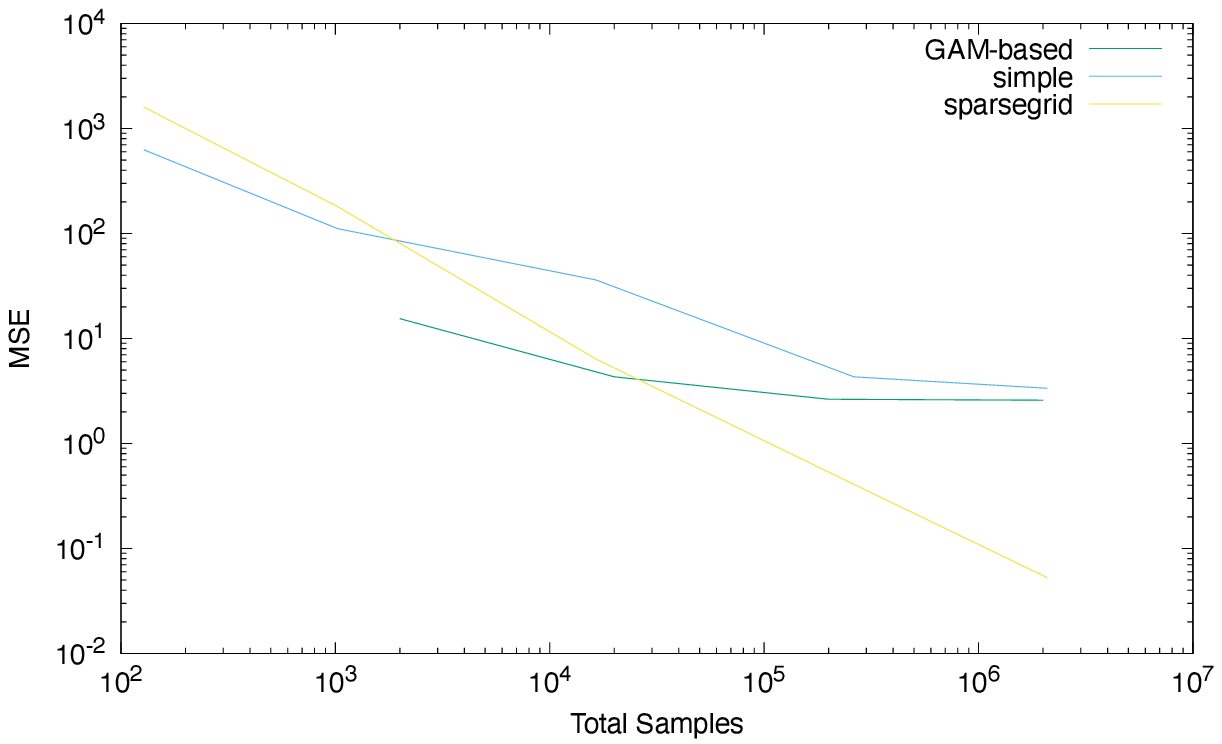}
         \caption{$n=1000$}
     \end{subfigure}
     \hfill
     \begin{subfigure}[b]{0.49\textwidth}
         \centering
         \includegraphics[width=\textwidth]{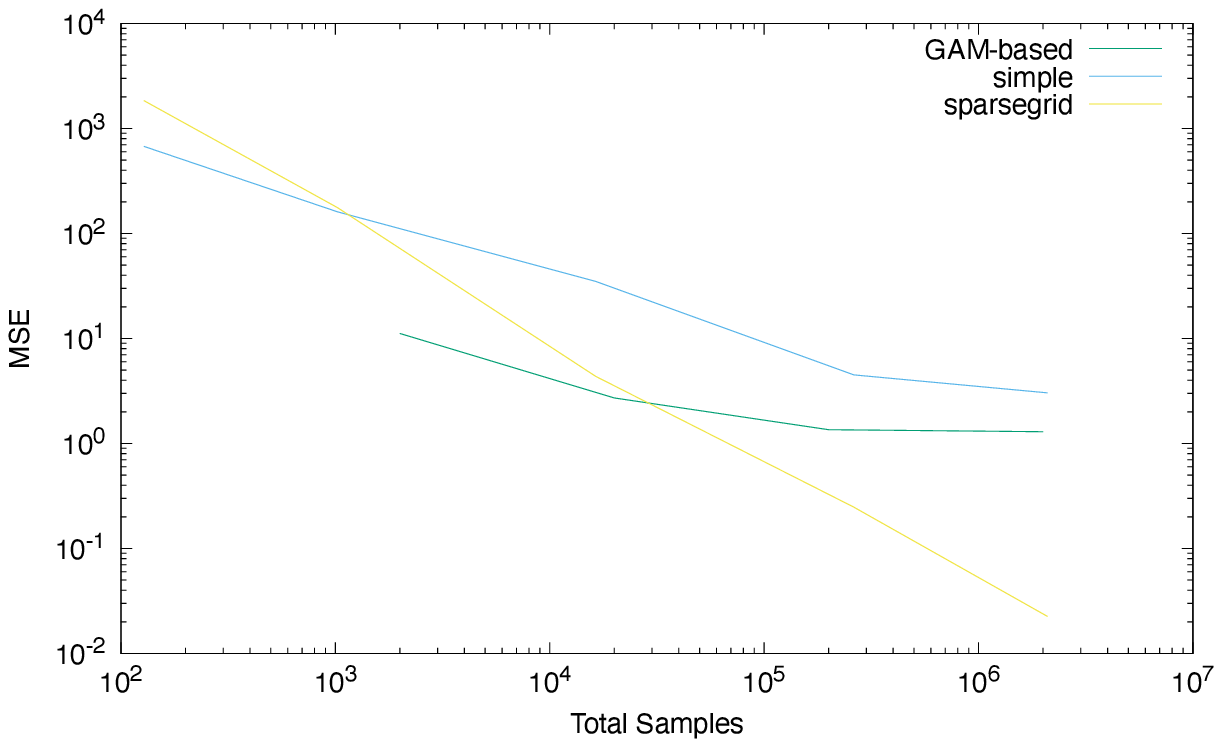}
         \caption{$n=2000$}
     \end{subfigure}
    \caption{Comparison of MSE between the three methods (our previous and proposed methods, and the GAM-based method) for the scenario EvSvGvA with various values of $n$. Our previous and proposed methods are denoted by {\tt simple} and {\tt sparsegrid}, respectively.}
    \label{fig:result_esga}
\end{figure}

The results for the EvSvGvA scenario with the same number of participants ($n=200, 500, 1000, 2000$) are shown in Figure~\ref{fig:result_esga}. In this scenario, the situation changes substantially; the GAM-based method exhibits inferior convergence behavior compared to our proposed method, and becomes even worse than our previous method for $n=200, 500$, when the number of total samples is large. It is possible that the MSE of the GAM-based method does not necessarily converge to 0. In contrast, the MSE of our proposed method converges to 0 at a rate similar to that in the EvSvG scenario.

%% file: 50-discussion.tex
Let us now discuss the numerical results obtained in these experiments. Firstly, we compare the methods in terms of consistency. As we have seen, all three methods appear to be consistent for the scenario EvSvG, while the GAM-based method does not necessarily converge to the exact EVSI value for the first example and the scenario EvSvGvA. GAM is a method of approximating a multivariate function by a sum of univariate functions. In the context of nested expectations, the inner conditional expectation $\mathbb{E}_{\rho(X|Y)}X$ is approximated simply by $\sum_{k=1}^{K}g_k(Y_k)$, where $g_1,\ldots,g_K$ are regressed functions. However, this additivity assumption does not hold in general, and thus an estimate of the GAM-based method may not converge to a correct value. In contrast, our proposed method is shown consistent under some conditions in the theoretical analysis, even if the additivity assumption is not valid.

Next, we compare the methods in terms of convergence. Our theoretical analysis demonstrates that our proposed method has an upper bound on the MSE of $O(N^{-1/K}(\log N)^2)$, while our previous method has an upper bound of $O(N^{-1/(2K)})$ under the same conditions. Thus, we expect our proposed method to converge faster than the previous method. In the numerical experiments, we observe that the MSE convergence slope for the proposed method is steeper than that for the previous method in all figures. Although the empirical slope does not agree with what is expected from Theorem~\ref{thm:main}, it is important to note that the theoretical analysis only provides an upper bound on the MSE, which may not be tight. Therefore, further research may be needed to refine the theoretical analysis and provide a more accurate estimate of the MSE convergence rate.

Finally, let us discuss the computational complexity of the methods, as the results presented so far are all in terms of MSE against the number of total samples $N$. As discussed in Remark~\ref{rem:dom}, the dominant parts of the process are assumed to be generating random samples from the joint distribution (of $\theta$ and $Y$ in our test problem) and evaluating the function $f$ (i.e., $\mathrm{NB}_d$ for all $d\in D$ in our test problem). Therefore, we expect the difference in computational complexity between the proposed method and other methods to be slight, as all methods used in our experiments have a quasilinear computational complexity in terms of $N$. In fact, our numerical experiments did not reveal any significant difference in computational time between the methods. Regarding the regression-based methods, it is possible to improve accuracy by adding interaction terms to the generalized additive models \cite{lou2013accurate} or using Gaussian process regression. However, there is a trade-off between accuracy improvement and computational cost, which makes it challenging to choose an appropriate model. In particular, Gaussian process regression has a computational cost of $O(N^3)$, limiting its applicability in regimes where $N$ is large.

%% file: 60-conclusion.tex
\section{Concluding Remarks}\label{sec:conclusion}

In this paper, we proposed a novel Monte Carlo estimator for nested expectations that does not require sampling from the inner conditional distribution. Our proposed method builds upon the previous method by the authors and is inspired by sparse grid quadrature. The theoretical analysis showed an upper bound on the MSE of $O(N^{-1/K}(\log N)^2)$ of our proposed method under mild conditions. Through a series of numerical experiments related to value of information analysis, we demonstrated the effectiveness of our method even with moderate dimensionality of the outer random variables.

Furthermore, we compared our proposed method with the previous method and the regression-based method in terms of convergence and found that our method may attain faster convergence. We also discussed the computational complexity of the methods and found that the difference in computational cost between the proposed method and the other methods is not significant, given that the dominant computational parts of the process are sampling from the joint distribution and evaluating the function values.

In conclusion, our proposed method provides a promising approach for efficiently estimating nested expectations with moderate dimensionality, as it is more efficient than the existing methods in terms of convergence and does not require sampling from the inner conditional distribution. Future research could investigate the theoretical properties in greater detail, as we have found some gaps between the empirical and theoretical convergence rates, and explore extensions of our method to handle more high-dimensional problems.

%% file: 90-backmatter.tex
\section*{Acknowledgements}
The work of the second author is supported by JSPS KAKENHI Grant Number 23K03210.

\bibliographystyle{plain}
\bibliography{91-references}